\documentclass{amsart}
\usepackage{amssymb}

\usepackage{amscd}
\usepackage{thmdefs}

\input{tcilatex}
\theoremstyle{definition}
\theoremstyle{remark}
\numberwithin{equation}{section}

\input tcilatex

\begin{document}
\title[Direct Producted Noncommutative Probability Spaces ]{Free Probability on a Direct Product of Noncommutative Probability Spaces}
\author{Ilwoo Cho}
\address{}
\email{}
\thanks{}
\date{}
\dedicatory{}
\thanks{}
\maketitle

\begin{abstract}
In this paper, we will consider the free probability on a direct product of
noncommutative probability spaces. Let $I$ be a finite set and let $\{(A_{i},
$ $\varphi _{i}):$ $i$ $\in $ $I\}$ be a collection of noncommutative
probability spaces, where $A_{i}$'s are unital algebras and $\varphi _{i}$'s
are linear functionals, for $i$ $\in $ $I.$ Define the direct product $A$ $=$
$\underset{i\in I}{\times }\ A_{i}$ of algebras $A_{i}$ and define the
conditional expectation $E$ $=$ $\underset{i\in I}{\times }$ $\varphi _{i}$
from $A$ onto the subalgebra $D_{\left| I\right| },$ generated by all
diagonal matrices in the matricial algebra $M_{\left| I\right| }(\Bbb{C}).$
We will consider the noncommutatibe probability space $(A,$ $E)$ with
amalgamation over $D_{\left| I\right| }.$ We observe that we can extend
almost all free probabilistic results in the scalar-valued case to those in
the $D_{\left| I\right| }$-valued case.
\end{abstract}

\strut \strut

Free Probability has been developed by various mathematician from 1980's.
There are two approaches to study it. One of them is the original
Voiculescu's pure analytic approach (See [3]) and the other one is the
Speicher and Nica's combinatorial approach (See [1], [10] and [11]). We will
use the Speicher's combinatorial amalgamated Free Probability introduced and
studied in [11]. Let $B$ $\subset $ $A$ be unital algebras with $1_{A}=1_{B}$
(equivalently, $A$ is an algebra over $B$). Suppose that there exists a
conditional expectation $E$ $:$ $A$ $\rightarrow $ $B$ satisfying the
bimodule map property and

\strut

(i) \ \ $E(b)=b,$ for all $b\in B$

(ii) \ $E(bab^{\prime })=bE(a)b^{\prime },$ for all $b,b^{\prime }\in B$ and 
$a\in A.$

\strut

Then the algebraic pair $(A,E)$ is called a noncommutative probability space
with amalgamation over $B$ (or an amalgamated noncommutative probability
space over $B.$ See [11]). All elements in $(A,$ $E)$ are said to be $B$%
-valued random variables. When $B=\Bbb{C}$ and $E$ is a linear functional,
then we call this structure a (scalar-valued) noncommutative probability
space and the elements (free) random variables. Let $a$ $\in $ $(A,$ $E)$ be
a $B$-valued random variable. Then it contains the following (equivalent)
free probabilistic data,

\strut

\begin{center}
$E(b_{1}a...b_{n}a)$
\end{center}

and

\begin{center}
$k_{n}(b_{1}a,...,b_{n}a)\overset{def}{=}\underset{\pi \in NC(n)}{\sum }%
E_{\pi }(b_{1}a,...,b_{n}a)\mu (\pi ,1_{n}),$
\end{center}

\strut

called the $B$-valued $n$-th moment of $a$ and the $B$-valued $n$-th
cumulant of $a,$ respectively, for all $n$ $\in $ $\Bbb{N}$ and $b_{1},$
..., $b_{n}$ $\in B$ arbitrary, where $E_{\pi }$ $(...)$ is the
partition-dependent $B$-valued moment of $a$ and $NC(n)$ is the collection
of all noncrossing partitions and $\mu $ is the M\"{o}bius functional in the
incidence algebra $I_{2},$ as the convolution inverse of the zeta functional 
$\zeta $ defined by

\strut

\begin{center}
$\zeta (\pi _{1},\pi _{2})\overset{def}{=}\left\{ 
\begin{array}{lll}
1 &  & \text{if }\pi _{1}\leq \pi _{2} \\ 
0 &  & \text{otherwise,}
\end{array}
\right. $
\end{center}

\strut

for all $\pi _{1},$ $\pi _{2}\in NC(n),$ and $n\in \Bbb{N}$ (See [11]). When 
$b_{1}$ $=$ ... $=$ $b_{n}$ $=0_{B},$ for $n\in \Bbb{N},$ we say that $%
E(a^{n})$ and $k_{n}(a,$ $...,$ $a)$ are trivial $n$-th $B$-valued moment
and cumulant of $a,$ respectively. Recall that the collection $NC(n)$ is a
lattice with the total ordering $\leq $ (See [1], [2], [10] and [11]).

\strut

Let $A_{1}$ and $A_{2}$ be subalgebras of $A.$ We say that they are free
over $B$ if all mixed cumulants of $A_{1}$ and $A_{2}$ vanish. Let $S_{1}$
and $S_{2}$ be subsets of $A.$ We say that subsets $S_{1}$ and $S_{2}$ are
free over $B$ if the subalgebras $A_{1}$ $=$ $A\lg $ $\{S_{1},$ $B\}$ and $%
A_{2}$ $=$ $A\lg $ $\{S_{2},$ $B\}$ are free over $B.$ In particular, the $B$%
-valued random variables $a_{1}$ and $a_{2}$ are free over $B$ if the
subsets $\{a_{1}\}$ and $\{a_{2}\}$ are free over $B$ in $(A,$ $E).$
Equivalently, given two $B$-valued random variables $a_{1}$ and $a_{2}$ are
free over $B$ in $(A,$ $E)$ if all mixed $B$-valued cumulants of $a_{1}$ and 
$a_{2}$ vanish. (Recall that, when $A$ is a $*$-algebra, they are free over $%
B$ if all mixed $B$-valued cumulants of $P(a_{1},$ $a_{1}^{*})$ and $%
Q(a_{2}, $ $a_{2}^{*})$ vanish, for all $P,$ $Q$ in $\Bbb{C}$ $[z_{1},$ $%
z_{2}].$) \strut

\strut \strut

In this paper, we will consider Free Probability on a direct product of
finite noncommutative probability spaces, as a new model for noncommutative
probability space with amalgamation over the matricial subalgebra generated
by all diagonal matrices. The construction of them is very similar to those
of Toeplitz noncommutative probability spaces over Toeplitz matricial
algebras in [5]. Let $N$ $\in $ $\Bbb{N}$ and let

\strut

\begin{center}
$\mathcal{F}=\{(A_{i},\varphi _{i}):i=1,...,N\}.$
\end{center}

\strut

be the finite family of noncommutative probability spaces $(A_{1},$ $\varphi
_{1}),$ ..., $(A_{N},$ $\varphi _{N}),$ where $A_{1},$ ..., $A_{N}$ are
unital algebras and $\varphi _{1},$ ..., $\varphi _{N}$ are linear
functionals. Define the direct product $A$ of $A_{1},$ ..., $A_{N}$ by

\strut

\begin{center}
$A=\times _{j=1}^{N}A_{j}=\{(a_{1},...,a_{N}):a_{j}\in A_{j},$ $j=1,...,N\}.$
\end{center}

\strut

Now, define the diagonal matricial algebra $D_{N}$ in the matricial algebra $%
M_{N}(\Bbb{C})$ by

\strut

\begin{center}
$D_{N}=\Bbb{C}[\{e_{11},...,e_{NN}\}],$
\end{center}

\strut

where $\{e_{ij}:i,$ $j=1,...,N\}$ is the canonical matrix units of $M_{N}(%
\Bbb{C}).$ i.e.,

\strut

\begin{center}
$e_{ii}= 
\begin{array}{ll}
\text{ \ \ \ \ \ \ \ \ \ \ \ \ \ }i\text{-th} &  \\ 
\left( 
\begin{array}{lllll}
0 & \cdots & \downarrow & \cdots & 0 \\ 
& \ddots &  &  &  \\ 
\vdots &  & \frame{1} &  & \leftarrow \\ 
&  &  & \ddots &  \\ 
0 & \cdots &  & \cdots & 0
\end{array}
\right) & i\text{-th}
\end{array}
$
\end{center}

\strut

for all \thinspace $i=1,...,N.$

\strut

Then we can define the conditional expectation $E$ $:$ $A$ $\rightarrow $ $%
D_{N}$ by

\strut

\begin{center}
$E\left( (a_{1},...,a_{N})\right) =\left( 
\begin{array}{lll}
\varphi _{1}(a_{1}) &  & \text{ \ }0 \\ 
& \ddots &  \\ 
\text{ \ \ }0 &  & \varphi _{N}(a_{N})
\end{array}
\right) .$
\end{center}

\strut

Then the algebraic pair $\left( A,E\right) $ is a noncommutative probability
space with amalgamation over $D_{N}.$ For convenience, sometimes, we will
denote $E$ by $(\varphi _{1},$ ..., $\varphi _{N}).$ We will call it the
direct producted noncommutative probability space. Notice that $D_{N}$ is a
commutative subalgebra of $M_{N}(\Bbb{C}).$ The main purpose of this paper
is to study Free Probability on this direct producted noncommutative
probability space $(A,E).$

\strut

In Chapter 1, we will observe the free structure of the direct producted
noncommutative probability space $(A,E)$, for instance, the $D_{N}$-freeness
on $(A,$ $E).$ In particular, we can see that if $x$ $=$ $(a_{1},$ ..., $%
a_{N})$ is in $(A,$ $E),$ then

\strut

\begin{center}
$E(x^{n})=\left( 
\begin{array}{lll}
\varphi _{1}(a_{1}^{n}) &  & \text{ \ \ \ \ \ \ \ }0 \\ 
& \ddots &  \\ 
0 &  & \varphi _{N}(a_{N}^{n})
\end{array}
\right) ,$
\end{center}

and

\begin{center}
$k_{n}\left( \underset{n\text{-times}}{\underbrace{x,.....,x}}\right)
=\left( 
\begin{array}{lll}
k_{n}^{(1)}\left( \underset{n\text{-times}}{\underbrace{a_{1},...,a_{1}}}%
\right) &  & \text{ \ \ \ \ \ \ \ \ \ \ \ \ \ \ \ \ \ \ }0 \\ 
& \ddots &  \\ 
0 &  & k_{n}^{(N)}\left( \underset{n\text{-times}}{\underbrace{%
a_{N},...,a_{N}}}\right)
\end{array}
\right) ,$
\end{center}

\strut

for all $n\in \Bbb{N},$ where $k_{n}(...)$ is the $D_{N}$-valued cumulant
bimodule map with respect to the conditional expectation $E$ and $%
k_{n}^{(i)}(...)$ are scalar-valued cumulant functional with respect to the
linear functionals $\varphi _{i},$ for $i$ $=$ $1,$ ..., $N.$ The above $%
D_{N}$-valued cumulant relation shows that if $D_{N}$-valued random
variables $x$ $=$ $(a_{1},$ ..., $a_{N})$ and $y$ $=$ $(a_{1}^{\prime },$
..., $a_{N}^{\prime })$ are free over $D_{N}$ if and only if $a_{j}$ and $%
a_{j}^{\prime }$ are free in $(A_{j},$ $\varphi _{j}),$ for all $j$ $=$ $1,$
..., $N.$\strut

\strut \strut

In Chapter 2, we will consider the $D_{N}$-valued moment series and
R-transform of a $D_{N}$-valued random variable in $(A,$ $E).$ Since $D_{N}$
commutes with the direct product $A$ $=$ $\times _{j=1}^{N}$ $A_{j},$ we
only need to consider trivial $D_{N}$-valued cumulants of $D_{N}$-valued
random variables for studying the free probabilistic information (e.g., $%
D_{N}$-valued distribution etc) of them. So, as a formal series in $D_{N},$
we can define the R-transform $R_{x}(z)$ of a $D_{N}$-valued random variable 
$x$ $\in $ $(A,$ $E),$ by

\strut

\begin{center}
$R_{x}(z)=\sum_{n=1}^{\infty }k_{n}(x,...,x)\,z^{n}$ \ in \ $D_{N}[[z]],$
\end{center}

\strut

where $A[[z]]$ is a ring of formal series in $A$ and $z$ is an
indeterminent. Then, similar to the scalar-valued case (See [1] and [10]),
we have the following $D_{N}$-valued R-transform calculus;

\strut

(1) \ $R_{x+y}(z)=R_{x}(z)+R_{y}(z),$ in $D_{N}[[z]]$

(2) \ $R_{x,\,y}(z_{1},\,z_{2})=R_{x}(z_{1})+R_{y}(z_{2})$ in $%
D_{N}[[z_{1},z_{2}]]$

(3) \ $R_{xy}(z)=\left( R_{x}\,\,\frame{*}_{D_{N}}\,\,R_{y}\right) (z)$ \ in
\ $D_{N}[[z]].$

\strut

In (3), the symbol \frame{*}$_{D_{N}}$ means the $D_{N}$-valued \textbf{%
restricted} boxed convolution on $\Theta _{D_{N}}$, where $\Theta _{D_{N}}$
is the subset of $D_{N}$ $[[z]]$ with zero $D_{N}$-constant terms.

\strut

In Chapter 3, we will consider certain $D_{N}$-valued random variables. To
do that, we assume each noncommutative probability space $(A_{j},$ $\varphi
_{j})$ is a $W^{*}$-probability space, for $j$ $=$ $1,$ ..., $N.$ Then, by
using the product topology on $A$ $=$ $\times _{j=1}^{N}$ $A_{j},$ the
direct product $A$ is also a $W^{*}$-algebra and, by the subspace topology, $%
D_{N}$ is also a $W^{*}$-algebra. So, we can have a $W^{*}$-probability
space $(A,$ $E)$ with amalgamation over $D_{N}.$ We observe the $D_{N}$%
-semicircularity, $D_{N}$-evenness, $D_{N}$-valued R-diagonality, $D_{N}$%
-circularity and $D_{N}$-valued infinitely divisibility. By [P\_$D_{N}$], we
will represent the $D_{N}$-valued property of the $D_{N}$-valued random
variable $x$ $=$ $(a_{1},$ ..., $a_{N})$ and by [P\_$\Bbb{C}$], we will
denote the scalar-valued property of $a_{j}$ $\in $ $(A_{j},$ $\varphi
_{j}), $ for $j$ $=$ $1,$ ..., $N.$ Then, we have that

\strut

\begin{quote}
\textbf{The }$D_{N}$\textbf{-valued random variable }$x$\textbf{\ has }[P\_$%
D_{N}$]\textbf{\ if and only if each nonzero }$a_{j}$\textbf{\ has }[P\_$%
\Bbb{C}$].
\end{quote}

\strut

For example, $x$ $=$ $(a_{1},$ $0,$ $a_{3})$ is $D_{3}$-semicircular if and
only if $a_{1}$ and $a_{3}$ are semicircular in $(A_{1},$ $\varphi _{1})$
and $(A_{3},$ $\varphi _{3}),$ respectively.

\strut \strut

Finally in Chapter 4, we will define the $D_{N}$-valued S-transforms of $%
D_{N}$-valued random variables. The operator-version of S-transforms are
merely known. But in this special $D_{N}$-valued case, we can define the
S-transforms and we can get the similar results like in the scalar-valued
case. In particular, if $x$ and $y$ are $D_{N}$-valued random variables in $%
(A,$ $E)$ and if $E(x)$ and $E(y)$ are invertible elements in $D_{N},$ under
the usual multiplication, then we can define the $D_{N}$-valued S-transforms 
$S_{x}(z)$ and $S_{y}(z)$ of $x$ and $y$ in $D_{N}$ $[[z]]$. And if $x$ and $%
y$ are free over $D_{N}$ in $(A,$ $E),$ then

\strut

\begin{center}
$S_{xy}(z)=\left( S_{x}(z)\right) \left( S_{y}(z)\right) .$
\end{center}

\strut

\strut

\strut

\section{Direct Producted Noncommutative Probability Spaces}

\strut

\strut

Throughout this chapter, let's fix $N\in \Bbb{N}$ and the collection $%
\mathcal{F}$ of (scalar-valued) noncommutative probability spaces,

\strut

\begin{center}
$\mathcal{F}=\{(A_{i},\varphi _{i}):i=1,...,N\}.$
\end{center}

\strut

Then, for the given unital algebras $A_{1},$ ..., $A_{N},$ we can define the
direct producted unital algebra

$\strut $

\begin{center}
$A$ $=$ $\times _{j=1}^{N}$ $A_{j}=\{(a_{1},...,a_{N}):a_{j}\in A_{j},$ $%
j=1,...,N\},$
\end{center}

\strut

as a set with its usual vector addition and the vector multiplication
defined componentwisely by

\strut

\begin{center}
$(a_{1},...,a_{N})\cdot (a_{1}^{\prime },...,a_{N}^{\prime
})=(a_{1}a_{1}^{\prime },...,a_{N}a_{N}^{\prime }),$
\end{center}

\strut

for all $(a_{1},...,a_{N}),$ $(a_{1}^{\prime },...,a_{N}^{\prime })$ $\in $ $%
A.$ Then the vector multiplication on $A$ is associative and hence $A$ is
again a unital algebra with its unity $(1,$ ..., $1).$

\strut

Let $M_{N}(\Bbb{C})$ be the matricial algebra generated by all $N$ $\times $ 
$N$-matrices with its canonical matrix units $\{e_{ij}$ $:$ $i,$ $j$ $=$ $1,$
..., $N\},$ where

\strut

\begin{center}
$e_{ij}= 
\begin{array}{ll}
\text{ \ \ \ \ \ \ \ }j\text{-th} &  \\ 
\left( 
\begin{array}{llllll}
0 & \downarrow &  & \cdots &  & \text{ \ \ }0 \\ 
& \frame{1} &  &  &  & \longleftarrow \\ 
&  &  &  &  &  \\ 
\text{ }\vdots &  &  & \ddots &  & \text{ \ \ \ }\vdots \\ 
&  &  &  &  &  \\ 
0 &  &  &  &  & \text{ \ \ }0
\end{array}
\right) & 
\begin{array}{l}
i\text{-th}
\end{array}
\end{array}
$
\end{center}

\strut

Define the subalgebra $D_{N}$ of $M_{N}(\Bbb{C})$ by

\strut

\begin{center}
$D_{N}\overset{def}{=}\Bbb{C}[\{e_{jj}:j=1,...,N\}].$
\end{center}

\strut

Then $D_{N}$ is a commutative subalgebra of $M_{N}(\Bbb{C})$ and it is
generated by all $N$ $\times $ $N$ diagonal matrices in $M_{N}(\Bbb{C}).$ We
will call this subalgebra $D_{N},$ the $N$-th diagonal algebra of $M_{N}(%
\Bbb{C}).$ Notice that we can regard $D_{N}$ as the algebra $\Bbb{C}^{N}$
with its usual vector addition and the following vector multiplication,

\strut

\begin{center}
$\left( \alpha _{1},...,\alpha _{N}\right) \cdot \left( \alpha _{1}^{\prime
},...,\alpha _{N}^{\prime }\right) =\left( \alpha _{1}\alpha _{1}^{\prime
},...,\alpha _{N}\alpha _{N}^{\prime }\right) .$
\end{center}

\strut

From now, if we mention the $N$-th diagonal algebra $D_{N},$ then it is
regarded as the algebra $\left( \Bbb{C}^{N},+,\cdot \right) ,$ where the
vector multiplication ($\cdot $) is defined as above. Also, notice that,
under this assumption, the $N$-th diagonal algebra $D_{N}$ is a subalgebra
of the direct product $A$ $=$ $\times _{j=1}^{N}$ $A_{j}.$ Moreover,

\strut

\begin{center}
$1_{D_{N}}=(1,...,1)=1_{A}.$
\end{center}

\strut

Thus the direct product $A$ of $A_{1},$ ..., $A_{N}$ is an algebra over $%
D_{N}.$

\strut

\begin{definition}
The algebra $\left( \Bbb{C}^{N},+,\cdot \right) ,$ defined in the above
paragraph, is called the $N$-th diagonal algebra and we will denote it again
by $D_{N}.$ Define the direct producted noncommutative probability space $A$
of $\mathcal{F}$ by the noncommutative probability space $(A,$ $E)$ with
amalgamation over the $N$-th diagonal algebra $D_{N},$ where $A$ $=$ $\times
_{j=1}^{N}A_{j}$ is the direct product of $A_{1},$ ..., $A_{N}$ and $E$ $:$ $%
A$ $\rightarrow $ $D_{N}$ is the conditional expectation from $A$ onto $D_{N}
$ defined by

\strut 

$\ \ \ \ \ \ \ E\left( (a_{1},...,a_{N})\right) =\left( \varphi
_{1}(a_{1}),...,\varphi (a_{N})\right) ,$

\strut 

for all $(a_{1},...,a_{N})\in A.$ Sometimes, we will denote $E$ by $\left(
\varphi _{1},...,\varphi _{N}\right) .$
\end{definition}

\strut

It is easy to see that the $\Bbb{C}$-linear map $E$ is indeed a conditional
expectation;

\strut

(i) \ \ $E\left( (\alpha _{1},...,\alpha _{N})\right) =\left( \varphi
_{1}(\alpha _{1}),...,\varphi _{N}(\alpha _{N})\right) =\left( \alpha
_{1},...,\alpha _{N}\right) ,$

\strut

\ \ \ \ \ \ for all $(\alpha _{1},...,\alpha _{N})\in D_{N}.$

\strut

(ii) \ $E\left( (\alpha _{1},...,\alpha _{N})\left( a_{1},...,a_{N}\right)
(\alpha _{1}^{\prime },...,\alpha _{N}^{\prime })\right) $

\strut

$\ \ \ \ \ \ \ =E\left( (\alpha _{1}a_{1}\alpha _{1}^{\prime },...,\alpha
_{N}a_{N}\alpha _{N}^{\prime })\right) =\left( \varphi (\alpha
_{1}a_{1}\alpha _{1}^{\prime }),...,\varphi (\alpha _{N}a_{N}\alpha
_{N}^{\prime })\right) $

\strut

$\ \ \ \ \ \ \ =\left( \alpha _{1},...,\alpha _{N}\right) \cdot \left(
\varphi (a_{1}),...,\varphi (a_{N})\right) \cdot (\alpha _{1}^{\prime
},...,\alpha _{N}^{\prime })$

\strut

$\ \ \ \ \ \ \ =(\alpha _{1},...,\alpha _{N})\left(
E((a_{1},...,a_{N}))\right) (\alpha _{1}^{\prime },...,\alpha _{N}^{\prime
}),$

\strut

\ \ \ \ \ \ for all $(\alpha _{1},...,\alpha _{N}),$ $(\alpha _{1}^{\prime
},...,\alpha _{N}^{\prime })$ $\in $ $D_{N}$ and $(a_{1},...,a_{N})\in A.$

\strut

By (i) and (ii), the map $E$ is a conditional expectation from $A$ $=$ $%
\times _{j=1}^{N}$ $A_{j}$ onto $D_{N}.$ Thus the algebraic pair $(A,$ $E)$
is a noncommutative probability space with amalgamation over the $N$-th
diagonal algebra $D_{N}.$

\strut

Now, we will consider the $D_{N}$-freeness on the direct producted
noncommutative probability space $(A,$ $E)$. Notice that the $N$-th diagonal
algebra $D_{N}$ satisfies that

\strut

(1.1) $\ \ \ \ \ \ \ dx=xd,$ for all $d\in D_{N}$ and $x\in A$,

\strut

as a subalgebra of our direct product $A$ $=$ $\times _{j=1}^{N}$ $A_{j}.$
By (1.1) and the commutativity of $D_{N},$ we only need to consider the
trivial $D_{N}$-valued moments and cumulants of $D_{N}$-valued random
variables to study the free probabilistic information of them. i.e., we have
that, for any $a\in A,$

\strut

(1.2) $\ \ \ \ \ \ \ \ \ \ \ E(d_{1}a...d_{n}a)=\left( d_{1}...d_{n}\right)
E(a^{n})$

\strut

and

\strut

(1.3) \ \ \ \ \ \ \ $k_{n}\left( d_{1}a,...,d_{n}a\right) =\left(
d_{1}...d_{n}\right) k_{n}\left( a,...,a\right) ,$

\strut

for all $n\in \Bbb{N}$ and for any arbitrary $d_{1},...,d_{n}\in D_{N}.$ 

\strut

\begin{proposition}
Let $(A,E)$ be the direct producted noncommutative probability space with
amalgamation over the $N$-th diagonal algebra $D_{N},$ where $A$ $=$ $\times
_{j=1}^{N}$ $A_{j}$ and let $a_{1}$ and $a_{2}$ be $D_{N}$-valued random
variables in $(A,$ $E).$ Then they are free over $D_{N}$ in $(A,$ $E)$ if
and only if all mixed \textbf{trivial }$D_{N}$-valued cumulants of them
vanish. $\square $
\end{proposition}

\strut

We could find the similar fact in [5], for the Toeplitz noncommutative
probability spaces over the Toeplitz matricial algebras. We will consider
the $D_{N}$-valued moments of an arbitrary random variable in the direct
producted noncommutative probability space.

\strut

\begin{proposition}
Let $\left( A=\times _{j=1}^{N}A_{j},\text{ }E=(\varphi _{1},...,\varphi
_{N})\right) $ be the direct producted noncommutative probability space over
the $N$-th diagonal algebra $D_{N}$ and let $x$ $=$ $(a_{1},$ $...,$ $a_{N})$
be the $D_{N}$-valued random variable in $(A,$ $E).$ Then the trivial $n$-th
moment of $x$ is

\strut 

(1.4)$\ \ \ \ \ \ \ \ \ \ \ E(x^{n})=\left( \varphi
_{1}(a_{1}^{n}),...,\varphi _{N}(a_{N}^{n})\right) ,$

\strut 

for all $n\in \Bbb{N}.$ \ $\square $
\end{proposition}

\strut

The above proposition is easily proved by the straightforward computation,
for the fixed $n$ $\in $ $\Bbb{N}.$ It shows that if we know the $n$-th
moments of $a_{j}$ $\in $ $(A_{j},$ $\varphi _{j}),$ for $j$ $=$ $1,$ ..., $%
N,$ then we can compute the $D_{N}$-valued moment of the $D_{N}$-valued
random variable $(a_{1},$ ..., $a_{N})$ in $(A,$ $E).$ Now, let's compute
the trivial $n$-th cumulant of an arbitrary $D_{N}$-valued random variable;

\strut

\begin{proposition}
Let $\left( A=\times _{j=1}^{N}A_{j},E=(\varphi _{1},...,\varphi
_{N})\right) $ be the direct producted noncommutative probability space over
the $N$-th diagonal algebra $D_{N}$ and let $x$ $=$ $(a_{1},$ ..., $a_{N})$
be the $D_{N}$-valued random variable in $(A,$ $E).$ Then the trivial $n$-th
cumulant of $x$ is

\strut \strut 

\ \ $\ \ \ k_{n}\left( \underset{n\text{-times}}{\underbrace{x,........,x}}%
\right) =\left( k_{n}^{(1)}\left( a_{1},...,a_{1}\right)
,...,k_{n}^{(N)}(a_{N},...,a_{N})\right) ,$

\strut 

for all $n\in \Bbb{N},$ where $k_{n}^{(i)}(...)$ is the $n$-th cumulant
functional with respect to the noncommutative probability space $(A_{i},$ $%
\varphi _{i}),$ for all $i$ $=$ $1,$ ..., $N.$
\end{proposition}

\strut

\begin{proof}
Fix $n\in \Bbb{N}.$ Then

\strut

$\ \ k_{n}\left( x,...,x\right) =\underset{\pi \in NC(n)}{\sum }E_{\pi
}(x,...,x)\mu (\pi ,1_{n})$

\strut

(1.5)\ \ \ \ 

\strut

$\ \ \ =\underset{\pi \in NC(n)}{\sum }\left( \underset{V\in \pi }{\Pi }%
E_{V}(x,...,x)\right) \mu (\pi ,1_{n})$

\strut

by (1,1), where $E_{V}(x,...,x)=E(x^{\left| V\right| }),$ where $\left|
V\right| $ is the length of the block (See [1] and [11])

\strut

$\ \ \ =\underset{\pi \in NC(n)}{\sum }\left( \underset{V\in \pi }{\Pi }%
\left( (\varphi _{1}(a_{1}^{\left| V\right| }),...,\varphi
_{N}(a_{N}^{\left| V\right| }))\right) \right) \mu (\pi ,1_{n})$

\strut

$\ \ \ =\underset{\pi \in NC(n)}{\sum }\left( \left( \underset{V\in \pi }{%
\Pi }\varphi _{1}(a_{1}^{\left| V\right| }),...,\underset{V\in \pi }{\Pi }%
\varphi _{N}(a_{N}^{\left| V\right| })\right) \right) \mu (\pi ,1_{n})$

\strut

$\ \ \ =\underset{\pi \in NC(n)}{\sum }\left( \left( (\underset{V\in \pi }{%
\Pi }\varphi _{1}(a_{1}^{\left| V\right| }))\mu (\pi ,1_{n}),...,(\underset{%
V\in \pi }{\Pi }\varphi _{N}(a_{N}^{\left| V\right| }))\mu (\pi
,1_{n})\right) \right) ,$

\strut

since the $N$-th diagonal algebra $D_{N}$ is a vector space (i.e., if we let 
$\mu _{\pi }$ $=$ $\mu (\pi ,$ $1_{n})$ in $\Bbb{C},$ for the fixed $\pi $ $%
\in $ $NC(n),$ then $(a_{1},$ ..., $a_{N})$ $\cdot $ $\mu _{\pi }$ $=$ $%
(a_{1}\mu _{\pi },$ ..., $a_{N}\mu _{\pi }),$ for all $(a_{1},$ ..., $a_{N})$
$\in $ $A$ $=$ $\times _{j=1}^{N}$ $A_{j}.$)

\strut

$\ \ \ =\left( \underset{\pi \in NC(n)}{\sum }\left( \underset{V\in \pi }{%
\Pi }\varphi _{1}(a_{1}^{\left| V\right| })\right) \mu (\pi ,1_{n}),...,%
\underset{\pi \in NC(n)}{\sum }\left( \underset{V\in \pi }{\Pi }\varphi
_{N}(a_{N}^{\left| V\right| })\right) \mu (\pi ,1_{n})\right) $

\strut

since $D_{N}$ is a vector space (i.e., $(a_{1},$ ..., $a_{N})$ $+$ $%
(a_{1}^{\prime },$ ..., $a_{N}^{\prime })$ $=$ $(a_{1}$ $+$ $a_{1}^{\prime
}, $ ..., $a_{N}$ $+$ $a_{N}^{\prime })$ in $A.$)

\strut

\strut (1.6)

\strut

$\ \ \ =\left(
k_{n}^{(1)}(a_{1},...,a_{1}),...,k_{n}^{(N)}(a_{N},...,a_{N})\right) ,$

\strut

where $k_{n}^{(i)}(...)$ is the (scalar-valued) $n$-th cumulant functional
with respect to the (scalar-valued) noncommutative probability space $%
(A_{i}, $ $\varphi _{i}),$ for all $i$ $=$ $1,$ ..., $N.$
\end{proof}

\strut

Therefore, the equality (1.6) in the above proposition shows that the $n$-th 
$D_{N}$-valued cumulant

$\strut $

\begin{center}
$k_{n}\left( (a_{1},...,a_{N}),...,(a_{1},...,a_{N})\right) $
\end{center}

\strut

of the $D_{N}$-valued random variable $(a_{1},...,a_{N})$ in the direct
producted noncommutative probability space $(A,$ $E),$ is nothing but the $N$%
-tuple of $n$-th (scalar-valued) cumulants of $a_{1},$ ..., $a_{N},$

\strut

\begin{center}
$\left( k_{n}^{(1)}(a_{1},...,a_{1}),...,k_{n}^{(N)}(a_{N},...,a_{N})\right)
.$
\end{center}

\strut

By (1.6), we can get the following $D_{N}$-freeness characterization on the
direct producted noncommutative probability space

$\strut $

\begin{center}
$\left( A=\times _{j=1}^{N}A_{j},\text{ }E=(\varphi _{1},...,\varphi
_{N})\right) .$
\end{center}

\strut

\begin{theorem}
Let $\left( A=\times _{j=1}^{N}A_{j},E=(\varphi _{1},...,\varphi
_{N})\right) $ be the given direct producted noncommutative probability
space over the $N$-th diagonal algebra $D_{N}$ and let $x_{1}$ $=$ $(a_{1},$
..., $a_{N})$ and $x_{2}$ $=$ $(b_{1},$ ..., $b_{N})$ be the $D_{N}$-valued
random variables in $(A,E).$ Then $x_{1}$ and $x_{2}$ are free over $D_{N}$
in $(A,E)$ if and only if nonzero $a_{j}$ and $b_{j}$ are free in $(A_{j}$, $%
\varphi _{j}),$ for all $j$ $=$ $1,$ ..., $N.$
\end{theorem}

\strut

\begin{proof}
($\Leftarrow $) Assume that random variables $a_{j}$ and $b_{j}$ are free in 
$(A_{j},\varphi _{j}),$ for all $j$ $=$ $1,$ ..., $N.$ Then, by the
freeness, all mixed $n$-th cumulants of $a_{j}$ and $b_{j}$ vanish, for all $%
j$ $=$ $1,$ ..., $N$ and for all $n$ $\in $ $\Bbb{N}$ $\setminus $ $\{1\}.$
By the previous theorem, it is sufficient to show that the $N$-tuples $x_{1}$
$=$ $(a_{1},$ ..., $a_{N})$ and $x_{2}$ $=$ $(b_{1},$ ..., $b_{N})$ in $(A,$ 
$E)$ have vanishing mixed \textbf{trivial} $D_{N}$-valued cumulants. Fix $n$ 
$\in $ $\Bbb{N}$ $\setminus $ $\{1\}$ and let $(x_{i_{1}},$ ..., $x_{i_{n}})$
are mixed $n$-tuple of $x_{1}$ and $x_{2},$ where $(i_{1},$ ..., $i_{n})$ $%
\in $ $\{1,$ $2\}^{n}.$ Then

\strut

$\ k_{n}$\strut $\left( x_{i_{1}},...,x_{i_{n}}\right) =k_{n}\left(
(a_{1}^{i_{1}},...,a_{N}^{i_{1}}),...,(a_{1}^{i_{n}},...,a_{N}^{i_{n}})%
\right) $

\strut

$\ \ \ \ \ \ \ \ \ \ \ =\left(
k_{n}^{(1)}(a_{1}^{i_{1}},...,a_{1}^{i_{n}}),...,k_{n}^{(N)}(a_{N}^{i_{1}},...,a_{N}^{i_{n}})\right) 
$

\strut

by the previous proposition

\strut

$\ \ \ \ \ \ \ \ \ \ \ =(0,...,0),$

\strut

by the hypothesis.

\strut

($\Rightarrow $) Let's assume that the $D_{N}$-valued random variables $x_{1}
$ $=$ $(a_{1},$ ..., $a_{N})$ and $x_{2}$ $=$ $(b_{1},$ ..., $b_{N})$ are
free over $D_{N}$ in $(A,$ $E)$ and assume also that there exists $j$ in $%
\{1,$ ..., $N\}$ such that nonzero $a_{j}$ and $b_{j}$ are not free in $%
(A_{j},$ $\varphi _{j}).$ Now, fix $n$ $\in $ $\Bbb{N}$ $\setminus $ $\{1\}$
and consider the mixed trivial $D_{N}$-valued cumulants of $x_{1}$ and $x_{2}
$.

\strut

$\ k_{n}\left( x_{i_{1}},...,x_{i_{n}}\right) =k_{n}\left(
(a_{1}^{i_{1}},...,a_{N}^{i_{1}}),...,(a_{1}^{i_{n}},...,a_{N}^{i_{n}})%
\right) $

\strut

$\ \ \ =\left(
k_{n}^{(1)}(a_{1}^{i_{1}},...,a_{1}^{i_{n}}),...,k_{n}^{(j)}(a_{j}^{i_{1}},...,a_{j}^{i_{n}}),..,k_{n}^{(N)}(a_{N}^{i_{1}},...,a_{N}^{i_{n}})\right) 
$

\strut

$\ \ \ =(0,...,0,$ $\ \underset{j\text{-th}}{\frame{?}}$ $\ ,0,...,0)\neq
(0,...,0,...,0),$

\strut

in general. Therefore, there exists the nonvanishing $D_{N}$-valued cumulant
of $x_{1}$ and $x_{2}.$ This contradict our assumption that $D_{N}$-valued
random variables $x_{1}$ and $x_{2}$ are free over $D_{N}$ in $(A,$ $E).$
\end{proof}

\strut

The above theorem shows that the $D_{N}$-freeness of $(a_{1},$ ..., $a_{N})$
and $(b_{1},$ ..., $b_{N})$ in the direct producted noncommutative
probability space $(A,$ $E)$ is characterized by the (scalar-valued)
freeness of nonzero $a_{j}$ and $b_{j}$ in $(A_{j},$ $\varphi _{j}),$ for
all $j$ $=$ $1,$ ..., $N.$

\strut

\begin{corollary}
Let $e_{i}=(0,...,0,a_{i},0,...,0)$ and $e_{j}=(0,...,0,a_{j},0,...,0)$ in $%
(A,$ $E),$ where $a_{k}$ $\in $ $(A_{k},$ $\varphi _{k}),$ for $k$ $=$ $i,$ $%
j.$ Then $e_{i}$ and $e_{j}$ are free over $D_{N}$ in $(A,$ $E),$ whenever $i
$ $\neq $ $j.$ \ $\square $
\end{corollary}

\strut

\strut Define subalgebras $A_{1}^{\prime },$ $...,$ $A_{N}^{\prime }$ of the
direct product $A=\times _{j=1}^{N}A_{j}$ by

\strut

\begin{center}
$A_{j}^{\prime }=\{(0,...,0,a_{j},0,...,0):a_{j}\in (A_{j},\varphi _{j})\},$
\end{center}

\strut

for all $j$ $=1,...,N.$ Then $A_{j}^{\prime }$ is the embedding of $A_{j}$
in $A.$ By the previous corollary, we can easily get the following;

\strut 

\begin{corollary}
The algebras $A_{1},$ ..., $A_{N}$ are free over $D_{N}$ in the direct
producted noncommutative probability space $(A,$ $E).$ $\square $
\end{corollary}

\strut

By definition, the unital algebra $A_{j}$ is always free from $\Bbb{C}$ (See
[1], [4], [10] and [11]). So, we can extend the above result as follows;

\strut

\begin{proposition}
The subalgebras $A_{j}^{\symbol{94}}=\Bbb{C}\times ...\times \Bbb{C}\times
A_{j}\times \Bbb{C}\times ...\times \Bbb{C}$, for $j$ $=$ $1,$ ..., $N,$ are
free from each other over $D_{N}$ in the direct producted noncommutative
probability space $(A,$ $E).$ $\square $
\end{proposition}

\strut \strut

\strut

\strut \strut

\section{$D_{N}$-Valued R-transform Calculus}

\strut

\strut

In this chapter, based on the facts in Chapter 1, we will consider the
amalgamated R-transform calculus of the direct producted noncommutative
probability space. The main purpose of this chapter is to extend the Nica's
scalar-valued R-transform calculus to the $D_{N}$-valued R-transform
calculus. Throughout this chapter, let $N$ $\in $ $\Bbb{N}$ and $(A_{j},$ $%
\varphi _{j})$ be (scalar-valued) noncommutative probability spaces, for $j$ 
$=$ $1,$ ..., $N.$ Also, let $A$ $=$ $\times _{j=1}^{N}$ $A_{j}$ be the
direct product algebra of $A_{1},$ ..., $A_{N},$ over the $N$-th diagonal
algebra $D_{N}$ and let $E$ $=$ $(\varphi _{1}$, ..., $\varphi _{N})$ be the
directed conditional expectation from the direct producted algebra $A$ onto
the algebra $D_{N}.$ The noncommutative probability space $(A,$ $E)$ with
amalgamation over $D_{N}$ is the direct proucted noncommutative probability
space. Now, on $(A,$ $E),$ we will define the $D_{N}$-valued moment series
and the $D_{N}$-valued R-transform of an arbitrary $D_{N}$-valued random
variable, as elements in the formal series in $D_{N}.$ i.e., they are
defined by the elements in $D_{N}$ $[[z]],$ where $z$ is an indeterminent.

\strut

\begin{definition}
Let $(A,E)$ be the given direct producted noncommutative probability space
over the $N$-th diagonal algebra $D_{N}$ and let $a$ $\in $ $(A,$ $E)$ be
the $D_{N}$-valued random variable. Define the $D_{N}$-valued moment series
of $a$ and the $D_{N}$-valued R-transform of $a$ by

\strut 

$\ \ \ \ \ \ \ \ \ \ \ \ \ \ \ \ \ \ \ \ \ M_{a}(z)=\sum_{n=1}^{\infty
}E(a^{n})\,z^{n}$

and

$\ \ \ \ \ \ \ \ \ \ \ \ \ \ \ \ \ R_{a}(z)=\sum_{n=1}^{\infty
}k_{n}(a,...,a)z^{n}$

\strut 

in $D_{N}[[z]],$ where $z$ is an arbitrary indeterminent and $E(a^{n})$ and $%
k_{n}(a,...,a)$ are the $D_{N}$-valued trivial $n$-th moment of $a$ and the $%
D_{N}$-valued trivial $n$-th cumulant of $a,$ respectively, in $B.$
\end{definition}

\strut

\begin{remark}
The above definition is little bit different from the original Voiculescu's
definition and from\ the Speicher's combinatorial definition. Recall that
Speicher defined the $D_{N}$-valued\ moment series of $a$ in $(A,$ $E)$ by

\strut 

$\ \ \ \ \ \ \ \ \ \ \ \ \ \ \ \sum_{n=1}^{\infty
}E(d_{1}ad_{2}a...d_{n}a)\in D_{N},$

\strut 

for arbitrary $d_{1},$ ..., $d_{n}$ $\in D_{N},$ for each $n\in \Bbb{N}.$
Also, he defined the $D_{N}$-valued R-transform of $a$ in $(A,$ $E)$ by

\strut 

$\ \ \ \ \ \ \ \ \ \ \ \ \ \ \ \sum_{n=1}^{\infty }k_{n}(d_{1}^{\prime
}a,...,d_{n}^{\prime }a)\in D_{N},$

\strut 

for arbitrary $d_{1}^{\prime },...,d_{n}^{\prime }\in D_{N},$ for each $n\in 
\Bbb{N}.$ However, in (1.2) and (1.3), we observed that it suffices to
consider the trivial $D_{N}$-valued moments and cumulants of $a.$ i.e., the $%
D_{N}$-valued trivial moments and cumulants contains the full free
probabilistic information of the $D_{N}$-valued random variable $a$ in $(A,$ 
$E).$ Therefore, instead of observing the general moments and cumulants, we
only observe the trivial ones. Thus, we defined the $D_{N}$-valued moment
series and the $D_{N}$-valued R-transform as in the above definition. Of
course, we do not need to define them as elements in the $D_{N}$-formal
series in $D_{N}$ $[[z]],$ but it makes us have the similar results with the
scalar-valued case in [1] and [10].
\end{remark}

\strut

Remember that, by Chapter 1, we have

\strut

\begin{center}
$E\left( (a_{1},...,a_{N})^{n}\right) =\left( \varphi
_{1}(a_{1}^{n}),...,\varphi _{N}(a_{N}^{n})\right) ,$
\end{center}

and

\begin{center}
$k_{n}\left( (a_{1},...,a_{N}),...,(a_{1},...,a_{N})\right) =\left(
k_{n}^{(1)}(a_{1},...,a_{1}),...,k_{n}^{(N)}(a_{N},...,a_{N})\right) ,$
\end{center}

\strut

for all $n\in \Bbb{N},$ where $(a_{1},...,a_{N})\in (A,E).$ Therefore, we
can get that if $a$ $=$ $(a_{1},$ ..., $a_{N})$ in the direct producted
noncommutative probability space $(A,$ $E),$ then the $D_{N}$-valued moment
series and the $D_{N}$-valued R-transform of $a$ are

\strut

(2.1) \ \ \ \ \ \ $M_{a}(z)=\sum_{n=1}^{\infty }\left( \varphi
_{1}(a_{1}^{n}),...,\varphi _{N}(a_{N}^{n})\right) z^{n}$

\strut

and

\strut

(2.2) \ \ $R_{a}(z)=\sum_{n=1}^{\infty }\left(
k_{n}^{(1)}(a_{1},...,a_{1}),...,k_{n}^{(N)}(a_{N},...,a_{N})\right) z^{n}$

\strut

in $D_{N}[[z]],$ respectively. Indeed, we have the similar R-transform
calculus with the scalar-valued case;

\strut

\begin{proposition}
Let $x=(a_{1},...,a_{N})$ and $y$ $=$ $(b_{1},...,b_{N})$ be $D_{N}$-valued
random variables in $(A,$ $E)$ and assume that they are free over $D_{N}$
(equivalently, for all $j$ $=$ $1,$ ..., $N,$ the random variables $a_{j}$
and $b_{j}$ in $(A_{j},$ $\varphi _{j})$ are free). Then

\strut 

(1) $\ R_{x+y}(z)=R_{x}(z)+R_{y}(z)$ \ in \ $D_{N}[[z]].$

(2) \ $R_{x,y}(z_{1},z_{2})=R_{x}(z_{1})+R_{y}(z_{2})$ in $%
D_{N}[[z_{1},z_{2}]].$ $\square $
\end{proposition}

\strut

Let $B$ be an arbitrary unital algebra and $\varphi $, the linear functional
on $B$ and let $(B,$ $\varphi )$ be the corresponding noncommutative
probability space. Assume that $x$ and $y$ be random variables in $(B,$ $%
\varphi )$ and suppose that $x$ and $y$ are free in $(B,$ $\varphi ).$ We
have the random variable $xy$ in $(B,$ $\varphi )$ and Speicher and Nica
showed that

\strut

(2.3) $\ \ \ k_{n}^{(\varphi )}\left( \underset{n\text{-times}}{\underbrace{%
xy,...,xy}}\right) =\underset{\pi \in NC(n)}{\sum }\left( k_{\pi }^{(\varphi
)}\left( \underset{n\text{-times}}{\underbrace{x,...,x}}\right) \cdot
k_{Kr(\pi )}^{(\varphi )}\left( \underset{n\text{-times}}{\underbrace{y,...,y%
}}\right) \right) ,$

\strut

for all $n\in \Bbb{N},$ where $k_{n}^{(\varphi )}(...)$ is the $n$-th
cumulant with respect to the given linear functional $\varphi $ on $B$ and $%
k_{\theta }^{(\varphi )}(...)$ is the partition-dependent cumulant, for all
noncrossing partition $\theta $ $\in $ $NC(n),$ and the map $Kr$ $:$ $NC(n)$ 
$\rightarrow $ $NC(n)$ is the Kreweras complementation map on $NC(n)$ (See
[1] and [10]). \strut

\strut \strut

We want to see that, on the direct producted noncommutative probability
space $(A,$ $E),$ a certain (operator-valued) relation like (2.3) holds
true. First, let's observe the following lemma;

\strut

\begin{lemma}
Let $n\in \Bbb{N}$ and $x=(a_{1},...,a_{N})\in (A,E)$, a $D_{N}$-valued
random variable and let $\pi $ $\in $ $NC(n)$ be a noncrossing partition.
Let $k_{\pi }$ $(x,$ ..., $x)$ be the partition-dependent cumulant bimodule
map of $x.$ Then

\strut 

$\ \ \ \ \ \ k_{\pi }(x,...,x)=\left( k_{\pi
}^{(1)}(a_{1},...,a_{1}),...,k_{\pi }^{(N)}(a_{N},...,a_{N})\right) .$
\end{lemma}

\strut

\begin{proof}
For $\pi \in NC(n),$

\strut

$\ k_{\pi }\left( \underset{n\text{-times}}{\underbrace{x,......,x}}\right) =%
\underset{V\in \pi }{\Pi }\left( k_{\left| V\right| }\left( \underset{\left|
V\right| \text{-times}}{\underbrace{x,........,x}}\right) \right) ,$

\strut

by (1.1)

\strut

$\ \ \ \ =\underset{V\in \pi }{\Pi }\left( k_{\left| V\right| }^{(1)}\left(
a_{1},...,a_{1}\right) ,...,k_{\left| V\right|
}^{(N)}(a_{N},...,a_{N})\right) ,$

\strut

by (1.6)

\strut

$\ \ \ \ =\left( \underset{V\in \pi }{\Pi }\left( k_{\left| V\right|
}^{(1)}(a_{1},...,a_{1})\right) ,...,\underset{V\in \pi }{\Pi }\left(
k_{\left| V\right| }^{(N)}(a_{N},...,a_{N})\right) \right) $

\strut

in $D_{N}$

\strut

$\ \ \ \ =\left( k_{\pi }^{(1)}(a_{1},...,a_{1}),...,k_{\pi
}^{(N)}(a_{N},...,a_{N})\right) .$

\strut
\end{proof}

\strut \strut

\begin{remark}
If $\pi =1_{n}$ in $NC(n),$ by the previous lemma, we have that

\strut 

$\ \ \ \ k_{1_{n}}\left( x,...,x\right) =\left(
k_{1_{n}}^{(1)}(a_{1},...,a_{1}),...,k_{1_{n}}^{(N)}(a_{N},...,a_{N})\right)
.$

\strut 

This is same as the formula (1.6).
\end{remark}

\strut \strut \strut

Assume now that $x=(a_{1},...,a_{N})$ and $y=(b_{1},...,b_{N})$ are $D_{N}$%
-valued random variables in the direct producted noncommutative probability
space $(A,$ $E).$ Also, assume that the $D_{N}$-valued random variables $x$
and $y$ are free over $D_{N}$ in $(A,$ $E).$ i.e., the nonzero random
variables $a_{j}$ and $b_{j}$ are free in $(A_{j},$ $\varphi _{j}),$ for all 
$j$ $=$ $1,$ ..., $N.$ Consider the $D_{N}$-valued random variable $xy$ in $%
(A,$ $E)$,

\strut

\begin{center}
$xy=\left( a_{1}b_{1},...,a_{N}b_{N}\right) .$
\end{center}

\strut

For the fixed $n\in \Bbb{N},$ we have that

\strut

$\ k_{n}\left( xy,...,xy\right) $

\strut

$\ \ \ =k_{n}\left(
(a_{1}b_{1},...,a_{N}b_{N}),...,(a_{1}b_{1},...,a_{N}b_{N})\right) $

\strut

by the previous chapter

\strut

$\ \ \ =\left( k_{n}^{(1)}\left( \underset{n\text{-times}}{\underbrace{%
a_{1}b_{1},...,a_{1}b_{1}}}\right) ,...,k_{n}^{(N)}\left( \underset{n\text{%
-times}}{\underbrace{a_{N}b_{N},...,a_{N}b_{N}}}\right) \right) $

\strut (2.4)

$\ \ \ =($ $\ \underset{\pi _{1}\in NC(n)}{\sum }k_{\pi _{1}}^{(1)}\left( 
\underset{n\text{-times}}{\underbrace{a_{1},...,a_{1}}}\right) \cdot
k_{Kr(\pi _{1})}^{(1)}\left( \underset{n\text{-times}}{\underbrace{%
b_{1},...,b_{1}}}\right) ,$

\strut

$\ \ \ \ \ \ \ \ \ \ \ ...,\,\,\,\underset{\pi _{N}\in NC(n)}{\sum }k_{\pi
_{N}}^{(N)}\left( \underset{n\text{-times}}{\underbrace{a_{N},...,a_{N}}}%
\right) \cdot k_{Kr(\pi _{N})}^{(N)}\left( \underset{n\text{-times}}{%
\underbrace{b_{N},...,b_{N}}}\right) $ $\ ),$

\strut

by (2.3). Then, by the rearrangement of the summation for each $j$-th entry (%
$j$ $=$ $1,$ ..., $N$), the formula (2.4) is equal to

\strut

(2.5)

\strut

$\ \ \ \ \ \ \ ($ $\ \underset{\pi \in NC(n)}{\sum }k_{\pi }^{(1)}\left( 
\underset{n\text{-times}}{\underbrace{a_{1},...,a_{1}}}\right) \cdot
k_{Kr(\pi )}^{(1)}\left( \underset{n\text{-times}}{\underbrace{%
b_{1},...,b_{1}}}\right) ,$

\strut

\ \ \ $\ \ \ \ \ \ \ \ \ \ \ ...,\,\,\,\underset{\pi \in NC(n)}{\sum }k_{\pi
}^{(N)}\left( \underset{n\text{-times}}{\underbrace{a_{N},...,a_{N}}}\right)
\cdot k_{Kr(\pi )}^{(N)}\left( \underset{n\text{-times}}{\underbrace{%
b_{N},...,b_{N}}}\right) $ $\ )$

\strut

which is same as

\strut

$\ \ \underset{\pi \in NC(n)}{\sum }\ \ $ $\ \left( k_{\pi }^{(1)}\left( 
\underset{n\text{-times}}{\underbrace{a_{1},...,a_{1}}}\right) \cdot
k_{Kr(\pi )}^{(1)}\left( \underset{n\text{-times}}{\underbrace{%
b_{1},...,b_{1}}}\right) ,...\right. $

\strut

\ \ \ \ \ \ \ \ \ \ \ $\ \ \ \ \ \ \ \ \ \ \ \left. ...,\,\,\,k_{\pi
}^{(N)}\left( \underset{n\text{-times}}{\underbrace{a_{N},...,a_{N}}}\right)
\cdot k_{Kr(\pi )}^{(N)}\left( \underset{n\text{-times}}{\underbrace{%
b_{N},...,b_{N}}}\right) \right) $

\strut

in $D_{N}$

\strut

$\ =\underset{\pi \in NC(n)}{\sum }\ \ ($ $\left( \ k_{\pi }^{(1)}\left( 
\underset{n\text{-times}}{\underbrace{a_{1},...,a_{1}}}\right)
,...,\,\,\,k_{\pi }^{(N)}\left( \underset{n\text{-times}}{\underbrace{%
a_{N},...,a_{N}}}\right) \right) $

\strut

$\ \ \ \ \ \ \ \ \ \ \ \ \ \ \ \ \ \ \ \cdot $ $\ \left( k_{Kr(\pi
)}^{(1)}\left( \underset{n\text{-times}}{\underbrace{b_{1},...,b_{1}}}%
\right) ,...,k_{Kr(\pi )}^{(N)}\left( \underset{n\text{-times}}{\underbrace{%
b_{N},...,b_{N}}}\right) \right) )$

\strut

in $D_{N}$

\strut

$\ =\underset{\pi \in NC(n)}{\sum }\left( k_{\pi }(x,...,x)\right) \left(
k_{Kr(\pi )}(y,...,y)\right) ,$

\strut

by the previous lemma. Therefore, we can get the following theorem;

\strut

\begin{theorem}
\strut Let $x$ and $y$ be $D_{N}$-valued random variables in $(A,E)$ and
assume that they are free over $D_{N}$ in $(A,$ $E).$ Then we have that

\strut 

(2.6) $\ \ \ k_{n}\left( xy,...,xy\right) =\underset{\pi \in NC(n)}{\sum }%
\left( k_{\pi }(x,...,x)\right) \left( k_{Kr(\pi )}(y,...,y)\right) ,$

\strut 

for all $n\in \Bbb{N}$. $\square $
\end{theorem}

\strut \strut

The formula (2.6) is similar to the scalar-valued case (2.3). By (2.6), we
can get the following $D_{N}$-valued R-transform calculus;

\strut

\begin{proposition}
Let $x$ and $y$ be $D_{N}$-valued random variables in $(A,E)$ which are free
over $D_{N}.$ Then

\strut 

$\ \ \ \ \ \ \ \ \ \ \ \ \ \ \ R_{xy}(z)=\sum_{n=1}^{\infty }d_{n}\,z^{n}$ \
in $\ D_{N}[[z]]$

with

\strut 

$\ \ \ d_{n}=\underset{\pi \in NC(n)}{\sum }\left( k_{\pi }(x,...,x)\right)
\left( k_{Kr(\pi )}(y,...,y)\right) ,$ $\ \forall n\in \Bbb{N}.$

$\square $
\end{proposition}

\strut

Let $\Theta _{t}$ be the subset of the set $\Bbb{C}[[t]]$ of all formal
series with zero constant terms. i.e.,

\strut

\begin{center}
$\Theta _{t}=\{\sum_{n=1}^{\infty }\alpha _{n}t^{n}:\alpha _{n}\in \Bbb{C}%
,\,\,\,\forall n\in \Bbb{N}\}.$
\end{center}

\strut

In [1] and [10], Speicher and Nica defined the boxed convolution \frame{*}
on $\Theta _{t},$ as the binary operation on $\Theta _{t}$ defined by

\strut

\begin{center}
$(g_{1},g_{2})\longmapsto g_{1}\,\,\,\frame{*}\,\,\,g_{2}$
\end{center}

with

\begin{center}
$\gamma _{n}=\underset{\pi \in NC(n)}{\sum }\alpha _{\pi }\cdot \beta
_{Kr(\pi )},$
\end{center}

\strut

where $g_{1}(t)=\sum_{n=1}^{\infty }\alpha _{n}t^{n},$ \ $%
g_{2}(t)=\sum_{n=1}^{\infty }\beta _{n}t^{n}$ and

\strut

\begin{center}
$\left( g_{1}\,\,\frame{*}\,\,g_{2}\right) (t)=\sum_{n=1}^{\infty }\gamma
_{n}\,\,t^{n}$
\end{center}

\strut

in $\Theta _{t},$ where

\strut

\begin{center}
$\alpha _{\theta }=\underset{V\in \theta }{\Pi }\alpha _{\left| V\right| }$
\ \ \ and \ \ $\beta _{\theta }=\underset{B\in \theta }{\Pi }\beta _{\left|
B\right| },$
\end{center}

\strut

for all $\theta \in NC(n),$ for $n\in \Bbb{N}.$

\strut

\begin{definition}
Let $D_{N}$ be the $N$-th diagonal algebra and let $D_{N}[[z]]$ be the ring
of the formal series in $D_{N}.$ Also, let $\Theta _{D_{N}}$ be the subset
of $D_{N}$ $[[z]]$ consisting of all formal series in $D_{N}$ with zero $%
D_{N}$-constant terms. Define the restricted $D_{N}$-valued boxed
convolution \frame{*}$_{D_{N}}$ on $\Theta _{D_{N}}$ by

\strut 

$\ \ \ \ \ \ \ \ \ \ \ \ \ \ \ \ \ \left( g_{1}\,\,\frame{*}%
_{D_{N}}\,\,g_{2}\right) (z)=\sum_{n=1}^{\infty }d_{n}z^{n}$

with

$\ \ \ \ \ \ \ \ \ \ \ \ \ \ \ \ \ d_{n}=\underset{\pi \in NC(n)}{\sum }%
a_{\pi }\cdot b_{Kr(\pi )}\in D_{N},$

\strut 

where $g_{1}(z)=\sum_{n=1}^{\infty }a_{n}z^{n}$, \ $g_{2}(z)=\sum_{n=1}^{%
\infty }b_{n}z^{n}$ in $\Theta _{D_{N}},$ and where

\strut 

$\ \ \ \ \ \ \ \ \ \ \ \ \ \ \ a_{\theta }=\underset{V\in \theta }{\Pi }%
a_{\left| V\right| }$ \ \ and \ \ $b_{\theta }=\underset{B\in \theta }{\Pi }%
b_{\left| B\right| },$

\strut 

for all $\theta \in NC(n),$ $n\in \Bbb{N}.$
\end{definition}

\strut \strut \strut 

The previous proposition can be rewritten as follows, by the restricted $%
D_{N}$-valued boxed convolution \frame{*}$_{D_{N}}$;

\strut

\begin{corollary}
Let $x$ and $y$ be $D_{N}$-valued random variables in $(A,E)$ which are free
over $D_{N}.$ Then

\strut 

$\ \ \ \ \ \ \ \ \ \ \ \ \ \ \ R_{xy}(z)=\left( R_{x}\,\,\,\frame{*}%
_{D_{N}}\,R_{y}\right) (z)$ \ in \ $\Theta _{D_{N}}.$

$\square $
\end{corollary}

\strut

Thus, like the scalar-valued case, we can get the following R-transform
calculus on the direct producted noncommutative probability space $(A,$ $E)$
over the $N$-th diagonal algebra $D_{N}$ ; if $x$ and $y$ are free over $%
D_{N}$ in $(A,$ $E),$ then

\strut

(1) \ $R_{x+y}(z)=R_{x}(z)+R_{y}(z)$ \ \ \ \ in \ $D_{N}[[z]],$

\strut (2) \ $R_{x,\,y}(z_{1},z_{2})=R_{x}(z_{1})+R_{y}(z_{2})$ \ in \ $%
D_{N}[[z_{1},z_{2}]],$

\strut (3) \ $R_{xy}(z)=\left( R_{x}\,\,\frame{*}_{D_{N}}\,\,R_{y}\right)
(z) $ \ \ in \ $D_{N}[[z]].$

\strut

Let's define the following formal series

\strut

\begin{center}
$Zeta(z)=\sum_{n=1}^{\infty }\left( 1_{D_{N}}\right) z^{n}$
\end{center}

and

\begin{center}
$Mob(z)=\sum_{n=1}^{\infty }\left( \mu (0_{n},1_{n})\cdot 1_{D_{N}}\right)
z^{n},$
\end{center}

\strut

in $\Theta _{D_{N}},$ where $1_{D_{N}}=\left( \underset{N\text{-times}}{%
\underbrace{1,.......,1}}\right) $ in $D_{N}.$ Then we can easily get the
following rewritten M\"{o}bius inversion;

\strut

\begin{proposition}
Let $x$ be a $D_{N}$-valued random variable in a direct producted
noncommutative probability space $(A,$ $E)$ over $D_{N}.$ Then

\strut 

(1) \ \ \ $M_{x}(z)=\left( R_{x}\ \frame{*}_{D_{N}}\ Zeta\right) (z)$

(2) \ \ \ $R_{x}(z)=\left( M_{x}\text{ \ \frame{*}}_{D_{N}}\text{ \ }%
Mob\right) (z),$

\strut 

in $D_{N}[[z]].$ $\square $
\end{proposition}

\strut

The above proposition is nothing but the reformulated M\"{o}bius inversion
in terms of the restricted $D_{N}$-valued boxed convolution.

\strut \strut 

\strut 

\section{Random Variables in a Direct Producted $W^{*}$-Probability Spaces}

\strut

\strut

Throughout this chapter, we will fix $N\in \Bbb{N}$ and let $(A_{j},$ $%
\varphi _{j})$ be a $W^{*}$-probability space, for all \ $j$ $=$ $1,$ ..., $%
N,$ where $A_{j}$'s are von Neumann algebras and $\varphi _{j}$'s are states
on $A_{j}$ satisfying that $\varphi _{j}$ $(a^{*})$ $=$ $\overline{\varphi
_{j}(a)},$ in $\Bbb{C},$ for all $a$ $\in $ $A_{j},$ for $j$ $=$ $1,$ ..., $%
N.$ Recall that $\varphi _{j}$ is called a trace if $\varphi _{j}(ab)$ $=$ $%
\varphi _{j}(ba),$ for all $a,$ $b$ in $A_{j},$ for $j$ $=$ $1,$ ..., $N.$
In this case, the $W^{*}$-probability space $(A_{j},$ $\varphi _{j})$ is
said to be a tracial $W^{*}$-probability space. In this chapter, we will
consider the direct producted noncommutative probability space $(A,$ $E)$
over the $N$-th diagonal algebra $D_{N},$ where $A$ $=$ $\times _{j=1}^{N}$ $%
A_{j}$ is the direct product of von Neumann algebras $A_{1},$ ..., $A_{N}$
and $E$ $=$ $(\varphi _{1},$ ..., $\varphi _{N})$ is the conditional
expectation satisfying that

\strut

\begin{center}
$E\left( (a_{1},...,a_{N})^{*}\right) =E\left( (a_{1},...,a_{N})\right)
^{*}, $
\end{center}

\strut

for all $(a_{1},...,a_{N})\in (A,E).$ On $A$,

\strut

\begin{center}
$(a_{1},...,a_{N})^{*}=\left( a_{1}^{*},...,a_{N}^{*}\right) ,$
\end{center}

\strut

where $a_{j}^{*}$ is the adjoint of $a_{j}$ in $A_{j},$ for all $(a_{1},$ $%
...,$ $a_{N})$ $\in $ $A.$ Then $A$ is indeed a $*$-algebra and hence, under
the finite product topology, $A$ is a von Neumann algebra, too. Also we have
that

\strut

\begin{center}
$
\begin{array}{ll}
E\left( (a_{1},...,a_{N})^{*}\right) \strut & =E\left(
(a_{1}^{*},...,a_{N}^{*})\right) \\ 
& =\left( (\varphi _{1}(a_{1}^{*}),...,\varphi _{N}(a_{N}^{*}))\right) \\ 
& =\left( (\overline{\varphi _{1}(a_{1})},...,\overline{\varphi _{N}(a_{N})}%
)\right) \\ 
& =\left( (\varphi _{1}(a_{1}),...,\varphi _{N}(a_{N}))\right) ^{*} \\ 
& =E\left( (a_{1},...,a_{N})\right) ^{*}.
\end{array}
$
\end{center}

\strut

Clearly, we can regard the $N$-th diagonal algebra $D_{N},$ as a von Neumann
algebra under the subspace topology $\frak{T}_{A}$ $\cap $ $D_{N}$, where $%
\frak{T}_{A}$ is the product topology of the topologies $\frak{T}_{A_{1}},$
..., $\frak{T}_{A_{N}}$ of $A_{1},$ ..., $A_{N}.$ Thus we have the inclusion
of von Neumann algebras $D_{N}$ $\subset $ $A$ and we also have the
conditional expectation $E$ satisfying the above involution condition.
Notice that the continuity of $E$ is preserved by the continuity of $\varphi
_{j}$'s.

\strut

\begin{definition}
The algebraic pair $(A,$ $E),$ where $A$ and $E$ are given in the previous
paragraph, is called the direct producted $W^{*}$-probability space over the 
$N$-th diagonal algebra $D_{N}.$
\end{definition}

\strut

Let's consider certain $D_{N}$-valued random variables in a direct producted 
$W^{*}$-probability space. 

\strut 

\begin{definition}
(i) \ \ $\ x\in (A,E)$ is $D_{N}$-semicircular if it is self-adjoint and the
only nonvanishing trivial $D_{N}$-valued cumulant of $x$ is the second one.

(ii) \ $\ x\in (A,E)$ is $D_{N}$-even if all trivial $D_{N}$-valued moments
of $x$ vanish.

(iii) $\ x\in (A,E)$ is $D_{N}$-valued R-diagonal, if the only nonvanishing
mixed trivial $D_{N}$-valued cumulants of $x$ and $x^{*}$ are

\strut 

$\ \ \ \ \ \ \ \ \ k_{2n}\left( x,x^{*},...,x,x^{*}\right) $ \ \ or \ \ $%
k_{2n}(x^{*},x,...,x^{*},x),$

\strut 

for all $n\in \Bbb{N}.$
\end{definition}

\strut \strut \strut

\begin{proposition}
Let $x$ $=(a_{1},...,a_{N})\in (A,E)$ be a $D_{N}$-valued random variable
and if nonzero $a_{i}$ are even (or R-diagonal) in $(A_{i},$ $\varphi _{i})$%
, for all $i$ $\in $ $\{1,$ ..., $N\},$ then $x$ $\in $ $(A,$ $E)$ is $D_{N}$%
-even (resp. R-diagonal).
\end{proposition}

\strut

\begin{proof}
The $D_{N}$-valued evenness and R-diagonality are determined by the $D_{N}$%
-valued cumulant relations by definition. Also, by (1.6), we have that

\strut

$\ \ \ \ \ k_{n}\left( x,...,x\right) =\left(
k_{n}^{(1)}(a_{1},...,a_{1}),...,k_{n}^{(N)}(a_{N},...,a_{N})\right) ,$

\strut

for all $n\in \Bbb{N}$. Since nonzero $a_{i}$ $\in $ $(A_{i},$ $\varphi _{i})
$ are even, for all $i$ $\in $ $\{1,$ ..., $N\},$ $a_{i}$ are self-adjoint
and there exist nonvanishing even cumulants of $a_{i}$'s and hence the $D_{N}
$-valued random variable $x$ is self-adjoint in $A$ and all odd trivial $%
D_{N}$-valued cumulants of $x$ vanish and there exist nonvanishing even
trivial $D_{N}$-valued cumulants of $x$. Thus $x$ is $D_{N}$-even in $(A,$ $%
E).$

\strut

Similarly, by the R-diagonality of nonzero $a_{j}$'s,

\strut

$\ \ \ k_{2n}\left( x,x^{*},...,x,x^{*}\right) $

\strut

(3.1)

\strut

$\ \ \ \ \ \ =\left(
k_{2n}^{(1)}(a_{1},a_{1}^{*},...,a_{1},a_{1}^{*}),...,k_{2n}^{(N)}(a_{N},a_{N}^{*},...,a_{N},a_{N}^{*})\right) , 
$

\strut

for all $n\in \Bbb{N}$ and there exist $n\in \Bbb{N}$ such that (3.1) does
not vanish. The similar result holds true for $k_{2n}$ $(x^{*},$ $x,$ ..., $%
x^{*},$ $x),$ for $n$ $\in $ $\Bbb{N}.$ Let's assume that, under the same
condition, the $D_{N}$-valued random variable $x$ is not $D_{N}$-valued
R-diagonal. Equivalently, there exists $(u_{1},$ ..., $u_{n})$ $\in $ $\{1,$ 
$*\}^{n}$ such that $u_{k}$ $=$ $u_{k+1},$ for some $k$ $=$ $1,$ ..., $n$ $-$
$1,$ and

\strut

(3.2)$\ \ \ \ \ \ \ \ \ \ \ \ k_{n}\left(
x^{u_{1}},...,x^{u_{k}},x^{u_{k+1}},...,x_{n}\right) \neq 0_{D_{N}}.$

\strut 

By (1.6), we have

\strut

(3.3)$\ \ k_{n}(x^{u_{1}},...,x^{u_{n}})=\left(
k_{n}^{(1)}(a_{1}^{u_{1}},...,a_{1}^{u_{n}}),...,k_{n}^{(N)}(a_{N}^{u_{1}},...,a_{N}^{u_{n}})\right) . 
$

\strut

By the R-diagonality of all nonzero $a_{j}\in (A_{j},\varphi _{j}),$ for $j$ 
$=$ $1,$ ..., $N,$ the formula (3.3) vanish. However, this contradict our
assumption (3.2). Therefore, $x$ is $D_{N}$-valued R-diagonal in $(A,$ $E).$
\end{proof}

\strut

The above proposition shows that the evenness and the R-diagonality of
nonzero $a_{j}$ $\in $ $(A_{j},$ $\varphi _{j}),$ for some $j$ $=$ $1,$ ..., 
$N,$ guarantee the $D_{N}$-valued evenness and R-diagonality of $(a_{1},$
..., $a_{N})$ in the direct producted $W^{*}$-probability space $(A,$ $E)$
over $D_{N},$ respectively. Then how about the converse?

\strut

\begin{theorem}
Let $x=(a_{1},...,a_{N})\in (A,E)$ be a $D_{N}$-valued random variable. The $%
D_{N}$-valued random variable $x$ is $D_{N}$-even (or $D_{N}$-valued
R-diagonal) if and only if all nozero random variable $a_{j}$ are even in $%
(A_{j},$ $\varphi _{j}),$ for $j$ $=$ $1,$ ..., $N.$
\end{theorem}

\strut

\begin{proof}
By the previous proposition, we only need to show the necessary condition.
It is sufficient to prove the case when we have $x$ $=$ $(a_{1},$ ..., $%
a_{N})$ in $(A,$ $E)$ with $a_{j}$ $\neq $ $0,$ for all $j$ $=1,$ ..., $N.$

\strut

(1) Let $x$ be $D_{N}$-even. Then, by definition, all odd $D_{N}$-valued
moments of $x$ vanish. Recall that

\strut

$\ \ \ \ \ \ \ \ \ \ \ \ \ \ \ \ \ \ \ \ \ \ E\left( x^{n}\right) =\left(
\varphi _{1}(a_{1}^{n}),...,\varphi _{N}(a_{N}^{n})\right) ,$

\strut

for all $n\in \Bbb{N}.$ Assume that there exists $j\in \{1,...,N\}$ such
that $a_{j}$ is not an even element in $(A_{j},$ $\varphi _{j}).$ i.e.,
there exists at least one odd number $n_{o}$ in $\Bbb{N}$ such that $\varphi
_{j}(a_{j}^{n_{o}})$ $\neq $ $0$ in $\Bbb{C}.$ Then we can get that

\strut

$\ \ \ \ \ \ \ \ E(x^{n_{o}})=\left( \varphi _{1}(a_{1}^{n_{o}}),...,\varphi
_{j}(a_{j}^{n_{o}}),...,\varphi _{N}(a_{N}^{n_{o}})\right) \neq 0_{D_{N}},$

\strut

This contradict our assumption that $x$ is $D_{N}$-even.

\strut

\strut (2) Suppose that $x$ is $D_{N}$-valued R-diagonal and assume that
there exists $j$ in $\{1,$ ..., $N\}$ such that $a_{j}$ is not R-diagonal in 
$(A_{j},$ $\varphi _{j}).$ So, there exists an $n_{0}$-tuple $(u_{1},$ ..., $%
u_{n_{0}})$ in $\{1,$ $*\}^{n_{0}},$ for some $n_{0}$ $\in $ $\Bbb{N},$ such
that it is not alternating and

\strut

$\ \ \ \ \ \ \ \ \ \ \ \ \ \ \ \ \ \ \ \ \ k_{n}^{(j)}\left(
a_{j}^{u_{1}},...,a_{j}^{u_{n_{0}}}\right) \neq 0$ \ in \ $\Bbb{C}.$

\strut

Then

\strut

$\ \ \ \ \ \ k_{n_{0}}\left( x^{u_{1}},...,x^{u_{n_{0}}}\right) $

\strut

$\ \ \ \ \ \ \ =\left(
k_{n_{0}}^{(1)}(a_{1}^{u_{1}},...,a_{1}^{u_{n_{0}}}),...,k_{n_{0}}^{(j)}(a_{j}^{u_{1}},...,a_{j}^{u_{n_{0}}}),\right. 
$

\strut

$\ \ \ \ \ \ \ \ \ \ \ \ \ \ \ \ \ \ \ \ \ \ \ \ \ \ \ \ \ \ \left.
...,k_{n_{0}}^{(N)}(a_{N}^{u_{1}},...,a_{N}^{u_{n_{0}}})\right) $

\strut

$\ \ \ \ \ \ \ \neq 0_{D_{N}}.$

\strut

This contradict our assumption that $x$ is $D_{N}$-valued R-diagonal.\strut
\end{proof}

\strut \strut \strut 

Recall that, as embedded $W^{*}$-subalgebras, $A_{1},$ ..., $A_{N}$ are free
from each other over $D_{N}$ in the direct product $A$ $=$ $\times _{j=1}^{N}
$ $A_{j}.$ So, a $D_{N}$-even element $x$ $=$ $(a_{1},$ ..., $a_{N})$ is the 
$D_{N}$-free sum of nonzero $D_{N}$-even elements $e_{i}$ $=$ $(0,$ ..., $%
a_{i},$ ..., $0),$ for all $i$ $=$ $1,$ ..., $N.$ i.e.,

\strut

\begin{center}
$x=\sum_{j=1}^{N}e_{j}$ \ \ and \ \ $e_{j}$'s are free from each other, for $%
j$ $=1,...,N.$
\end{center}

\strut

We have the same result when we replace $D_{N}$-evenness to $D_{N}$-valued
R-diagonality.

\strut

\strut  We have just observed that the $D_{N}$-valued evenness and
R-diagonality of a $D_{N}$-valued random variable $\strut x$ $=$ $(a_{1},$
..., $a_{N})$ in $(A,$ $E)$ is totally characterized by the scalar-valued
evenness and R-diagonality of nonzero random variables $a_{j}$ $\in $ $%
(A_{j},$ $\varphi _{j}),$ for all $j$ $=$ $1,\,$..., $N.$ The $D_{N}$%
-semicircularity of $x$ $=$ $(a_{1},$ ..., $a_{N})$ is similarly
characterized by the semicircularity of nonzero $a_{1},$ ..., $a_{N}.$

\strut \strut \strut

\begin{theorem}
Let $a_{j}\in (A_{j},\varphi _{j})$ be random variables, for all $j$ $=$ $1,$
..., $N.$ Then the $D_{N}$-valued random variable $x$ $=$ $(a_{1},$ ..., $%
a_{N})$ in $(A,$ $E)$ is $D_{N}$-semicircular if and only if all nonzero $%
a_{j}$'s are semicircular in $(A_{j},$ $\varphi _{j}),$ for \ $j$ $=$ $1,$
..., $N.$
\end{theorem}

\strut

\begin{proof}
($\Leftarrow $) By (1.6), it is clear.

\strut 

($\Rightarrow $) It suffices to show the case when we have nonzero $a_{j}$%
's, for all \ $j$ $=$ $1,$ ..., $N.$ Assume that $a_{j}$ $\neq $ $0,$ in $%
(A_{j},$ $\varphi _{j})$, for all $j$ $=$ $1,$ ..., $N$, and suppose that
the $D_{N}$-valued random variable $x$ $=$ $(a_{1},$ ..., $a_{N})$ is $D_{N}$%
-semicircular in $(A,$ $E).$ Let's assume that there exists $j_{0}$ $\in $ $%
\{1,$ ..., $N\}$ such that $a_{j_{0}}$ is not semicircular in $(A_{j_{0}},$ $%
\varphi _{j_{0}}).$ Then there exists nonvanishing $n_{0}$-th cumulant $%
k_{n_{0}}^{(j_{0})}(a_{j_{0}},...,a_{j_{0}})$ of $a_{j_{0}},$ where $n_{0}$ $%
\neq $ $2$ in $\Bbb{N}.$ observe that

\strut

$\ \ \ k_{n_{0}}\left( x,...,x\right) $

\strut

$\ \ \ =\left(
k_{n_{0}}^{(1)}(a_{1},...,a_{1}),...,k_{n_{0}}^{(j_{0})}(a_{j_{0}},...,a_{j_{0}}),...,k_{n_{0}}^{(N)}(a_{N},...,a_{N})\right) 
$

\strut

$\ \ \ \neq (0,...,0,...,0)=0_{D_{N}}.$

\strut

This contradict our assumption that $x$ is $D_{N}$-semicircular.
\end{proof}

\strut

Also, by regarding the von Neumann algebras $A_{1},$ ..., $A_{N},$ as
embedded $W^{*}$-subalgebras in the direct product $A$ $=$ $\times
_{j=1}^{N} $ $A_{j}$, a $D_{N}$-semicircular element $x$ $=$ $(a_{1},$ ..., $%
a_{N})$ is the $D_{N}$-free sum of $D_{N}$-semicircular elements $e_{j}$ $=$ 
$(0,$ ..., $a_{j},$ ..., $0),$ for $j$ $=$ $1,$ ..., $N.$

\strut

\begin{definition}
Let $x$ $\in $ $(A,$ $E)$ be a $D_{N}$-valued random variable. We say that $x
$ is infinitely divisible if there exist $D_{N}$-valued random variables $%
x_{n,1},$ ..., $x_{n,n}$ in $(A,$ $E),$ for each $n$ $\in $ $\Bbb{N},$ such
that (i) they are free from each other over $D_{N}$ in $(A,$ $E)$ and (ii)
they are $D_{N}$-valued identically distributed. i.e.,

\strut 

$\ \ \ \ k_{m}(x,...,x)=\sum_{j=1}^{n}k_{m}\left( x_{n,j},...,x_{n,j}\right)
=nk_{m}\left( x_{n,j},...,x_{n,j}\right) .$
\end{definition}

\strut \strut \strut 

We can recognize that the $D_{N}$-valued infinitely divisibility of $D_{N}$%
-valued random variables is also defined by the $D_{N}$-cumulant relation.
Therefore, we have that;

\strut 

\begin{theorem}
The $D_{N}$-valued random variable $x$ $=$ $(a_{1},$ ..., $a_{N})$ in $(A,$ $%
E)$ is $D_{N}$-valued infinitely divisible if and only if all nonzero $a_{j}$
are infinitely divisible in $(A_{j},$ $\varphi _{j}),$ for $j$ $=$ $1,$ ..., 
$N.$ $\square $
\end{theorem}

\strut \strut 

So, all properties of a $D_{N}$-valued random variable $x$ depending on the $%
D_{N}$-valued cumulants are characterized by the properties of scalar-valued
random variables (which are the components of $x$) depending on the
scalar-valued cumulants.\strut 

\strut \strut 

\strut \strut

\section{$D_{N}$-Valued S-transforms}

\strut

\strut

\strut

In this chapter, we will consider the $D_{N}$-valued S-transform theory on
the direct producted noncommutative probability space $(A,$ $E)$ of
noncommutative probability spaces $\{$ $(A_{j},$ $\varphi _{j}),$ $j$ $=$ $1,
$ ..., $N$ $\},$ over the $N$-th diagonal algebra $D_{N},$ where $A$ $=$ $%
\times _{j=1}^{N}$ $A_{j}$ is the direct product of unital algebras $A_{1},$
..., $A_{N}$ and $E$ $=$ $(\varphi _{1},$ ..., $\varphi _{N})$ is the
conditional expectation from $A$ onto $D_{N}.$ The scalar-valued
S-transforms are introduced by Voiculescu to study the free probabilistic
data of the products of two free (scalar-valued) random variables. Nica
defined the S-transforms combinatorially.

\strut

Let $(B,$ $\varphi )$ be a noncommutative probability space with its linear
functional $\varphi $ $:$ $B$ $\rightarrow $ $\Bbb{C}$ and let $b$ $\in $ $%
(B,$ $\varphi )$ be a random variable satisfying that $\varphi (b)$ $\neq $ $%
0.$ Then define the S-transform of $b$ by

\strut \strut

(4.1) \ \ \ \ \ \ \ \ \ \ $s_{b}(t)=\frac{1+z}{z}\,\,m_{b}^{<-1>}(t)=\frac{1%
}{z}\,r_{b}^{<-1>}(t),$

\strut \strut

in $\Bbb{C}[[t]],$ where $m_{b}(t)$ and $r_{b}(t)$ are the moment series and
R-transform of $b$ in $\Theta _{t}$, and where $m_{b}^{<-1>}(t)$ and $%
r_{b}^{<-1>}(t)$ is the inverse of $m_{b}(t)$ and $r_{b}(t),$ with respect
to the composition on $\Bbb{C}[[t]].$ Notice that, by the assumption that $%
\varphi (b)$ $\neq $ $0$ and by (4.1), the S-transform $s_{b}(t)$ of $b$ is
not contained in $\Theta _{t}$, since $s_{b}(t)$ contains the constant term $%
\frac{1}{\varphi (b)}.$ He showed that if $b_{1}$ and $b_{2}$ are random
variables in a noncomutative probability space $(B,$ $\varphi )$ satisfying
that $\varphi (b_{1})$ $\neq $ $0$ $\neq $ $\varphi (b_{2})$ in $\Bbb{C},$
and if the random variables $b_{1}$ and $b_{2}$ are free in $(B,$ $\varphi
), $ then

\strut

(4.2) \ \ \ $\ \ \ \ \ \ \ \ \ \ \ s_{b_{1}b_{2}}(t)=\left(
s_{b_{1}}(t)\right) \left( s_{b_{2}}(t)\right) .$

\strut

We will extend the S-transform relation (4.2) to the $D_{N}$-valued case.

\strut

Define a subset $A_{-1}$ of $A$ by

\strut

\begin{center}
$A_{-1}=\{(a_{1},...,a_{N})\in A:\varphi _{j}(a_{j})\neq 0$ in $A_{j},$ $%
\forall j=1,...,N\}.$
\end{center}

Also, define a subset $D_{N}^{-1}$ of $D_{N}$ by

\strut

\begin{center}
$D_{N}^{-1}=\{(\alpha _{1},...,\alpha _{N})\in D_{N}:\alpha _{j}\neq 0,$ $%
\forall j=1,...,N\}.$
\end{center}

\strut

Let $x=(a_{1},...,a_{N})\in A_{-1}.$ Then

\strut

(4.3) \ \ \ $\ \ \ \ \ \ \ \ E\left( x\right) =\left( \varphi
_{1}(a_{1}),...,\varphi _{N}(a_{N})\right) \in D_{N}^{-1}\subset D_{N}$

\strut \strut \strut

Notice that the algebraic structure $\left( D_{N}^{-1},\text{ }\cdot \right) 
$ is a group, where ($\cdot $) is the usual vector multiplication on $D_{N}.$
Consider the set of formal series

$\strut $

\begin{center}
$D_{N}^{inv}[[z]]=\left\{ \sum_{n=0}^{\infty }d_{n}z^{n}\in D_{N}[[z]]: 
\begin{array}{l}
d_{0}\in D_{N}^{-1}, \\ 
d_{n}\in D_{N},\forall n\in \Bbb{N}
\end{array}
\right\} ,$
\end{center}

\strut

as a subset of $D_{N}[[z]]$, where $z$ is the indeterminent. By the very
definition, for all $f$ $\in $ $D_{N}^{inv}$ $[[z]],$ the $D_{N}$-constant
term of $f$ is always contained in $D_{N}^{-1}$ and hence every element $f$
in $D_{N}^{inv}$ $[[z]]$ has its multiplication inverse $f^{-1}.$ (Remark
the difference between the composition inverse $f^{<-1>}$ and the
multiplication inverse $f^{-1}$.) Thus it is easy to see that the algebraic
structure $(D_{N}^{inv}$ $[[z]],$ $\cdot )$ is also a group. \strut 

\strut \strut

Recall the definitions of $\Theta _{D_{N}}$ and the restricted $D_{N}$%
-valued boxed convolution \frame{*}$_{D_{N}}$ (See Chapter 2).

\strut

\begin{definition}
Define the subset of $\Theta _{D_{N}}^{inv}$ of $\Theta _{D_{N}}$ by

\strut 

$\ \ \ \ \ \ \ \ \ \Theta _{D_{N}}^{inv}=\{g\in \Theta _{D_{N}}:g$ is \frame{%
*}$_{D_{N}}$-invertible$\}.$
\end{definition}

\strut

Then $\left( \Theta _{D_{N}}^{inv},\text{ \ \frame{*}}_{D_{N}}\right) $ is a
group, too, lying in $D_{N}[[z]].$ Indeed, we can have that, if $g_{1},$ $%
g_{2}$ and $g_{3}$ are in $\left( \Theta _{D_{N}}^{inv},\frame{*}%
_{D_{N}}\right) ,$ then we can construct the corresponding R-transforms $%
R_{x_{1}}(z),$ $R_{x_{2}}(z)$ and $R_{x_{3}}(z)$ of $x_{1}$, $x_{2}$ and $%
x_{3}$ in some noncommutative probability spaces $(A,$ $E)$ over $D_{N}$, \
such that $D_{N}a$ $=$ $aD_{N},$ for all $a$ $\in $ $A,$ where $x_{1},$ $%
x_{2}$ and $x_{3}$ are free from each other over $D_{N}.$

\strut

\begin{remark}
The choice of the above unital algebra $A,$ which are commutes with $D_{N},$
is totally depending on our definition of the restricted $D_{N}$-valued
boxed convolution, \frame{*}$_{D_{N}}$ on $\Theta _{D_{N}}.$ So, if $g$ $\in 
$ $(\Theta _{D_{N}}^{inv},$ \ \frame{*}$_{D_{N}}),$ then, by \frame{*}$%
_{D_{N}},$ we can regard $g$ as the R-transform $R_{x}(z)$ of the $D_{N}$%
-valued random variable $x$ in a noncommutative probability space $(B,$ $%
E_{D_{N}}^{B})$ with amalgamation over $D_{N},$ where $D_{N}$ commutes with $%
B.$
\end{remark}

\strut

So, if $g_{1},$ $g_{2}$ and $g_{3}$ are in $\left( \Theta _{D_{N}}^{inv},%
\frame{*}_{D_{N}}\right) ,$ then

\strut

\begin{center}
$g_{1}\,\frame{*}_{D_{N}}\,\left( g_{2}\,\frame{*}_{D_{N}}g_{3}\right)
=R_{x_{1}x_{2}x_{3}}(z)=\left( g_{1}\,\frame{*}_{D_{N}}g_{2}\right) \,\frame{%
*}_{D_{N}}g_{3}.$
\end{center}

\strut

Thus, the $D_{N}$-valued restricted boxed convolution \frame{*}$_{D_{N}}$ is
associative on $\Theta _{D_{N}}^{inv}.$

\strut

(ii) \ The formal series $1_{D_{N}}\cdot z\in \Theta _{D_{N}}^{inv}$ is the 
\frame{*}$_{D_{N}}$-identity.

\strut

(iii) By definition, every $g\in \Theta _{D_{N}}^{inv}$ has its \frame{*}$%
_{D_{N}}$-inverse.

\strut

By (i), (ii) and (iii), $\left( \Theta _{D_{N}}^{inv},\,\,\frame{*}%
_{D_{N}}\right) $ is indeed a group. Suppose now that

$\strut $

\begin{center}
$g(z)=\sum_{n=1}^{\infty }d_{n}z^{n}$ \ \ and \ \ $g^{-1}(z)=\sum_{n=1}^{%
\infty }d_{n}^{\prime }z^{n},$
\end{center}

\strut

where $g^{-1}$ is the \frame{*}$_{D_{N}}$-inverse of $g$ in $\Theta
_{D_{N}}^{inv}.$ By definition, we have that

\strut

\begin{center}
$\left( g\,\,\,\frame{*}_{D_{N}}\,\,\,g^{-1}\right) (z)=\sum_{n=1}^{\infty
}p_{n}\,z^{n}=\left( g^{-1}\,\,\frame{*}_{D_{N}}\,\,g\right) (z)$
\end{center}

with

\begin{center}
$\underset{\pi \in NC(n)}{\sum }\left( d_{\pi }\right) \cdot (d_{Kr(\pi
)}^{\prime })=p_{n}=\underset{\theta \in NC(n)}{\sum }(d_{\theta }^{\prime
})(d_{Kr(\theta )}),$
\end{center}

satisfying

\begin{center}
$p_{1}=1_{D_{N}}$ \ \ and \ \ $p_{n}=0_{D_{N}},$ $\forall n\in \Bbb{N}%
\,\setminus \,\{1\}.$
\end{center}

\strut

Thus we have that\strut

\begin{center}
$d_{1}d_{1}^{\prime }=1_{D_{N}}=d_{1}^{\prime }d_{1}.$
\end{center}

\strut This shows that

$\strut $

(4.4) $\ \ \ \ \ \ \ \ \ \ \ \ \ \ \ \ \ \ \ d_{1},\,\,d_{1}^{\prime }\in
D_{N}^{-1}\subset D_{N}.$

\strut \strut

Conversely, if $d_{1}\in D_{N}^{-1}$ for $g(z)=\sum_{n=1}^{\infty }d_{n}z^{n}
$ in $\Theta _{D_{N}},$ then we can select $g^{-1}(z)$ $=$ $d_{1}^{-1}$ $+$ $%
\sum_{n=2}^{\infty }$ $\delta _{n}(d_{1})$ $z^{n}$ in $\Theta _{D_{N}}^{inv},
$ where $\delta _{n}(d_{1})$ are elements in $D_{N}$ depending on $d_{1},$
for all $n$ $\in $ $\Bbb{N}$ $\setminus $ $\{1\}.$

\strut

\begin{theorem}
The groups $\left( \Theta _{D_{N}}^{inv},\text{ \ \frame{*}}_{D_{N}}\right) $
and $\left( D_{N}^{inv}[[z]],\text{ }\cdot \right) $ are homomorphic groups. 
$\square $
\end{theorem}

\strut

The above theorem is proved by the little modification of the proof of the
Theorem 14.3 in [1]. Again, we have to remark that if we do not have the
condition (1.1), we cannot use the proof in [1] to prove the above theorem.
In Theorem 14.3 in [1], Nica showed that the groups $\left( \Theta
_{t}^{inv},\,\frame{*}\right) $ and $\left( \Bbb{C}[[t]]_{inv},\,\,\cdot
\right) $ are isomorphic, by the isomorphism $\mathcal{F}$ from $\Theta
_{t}^{inv}$ onto $\Bbb{C}[[t]]_{inv},$

\strut

\begin{center}
$\mathcal{F}\left( f(t)\right) =\frac{1}{t}f^{<-1>}(t),$ for all $f\in
\Theta _{t}.$
\end{center}

\strut

But we only show that $(\Theta _{D_{N}}^{inv},$ $\ \frame{*}_{D_{N}})$ and $%
(D_{N}^{inv}$ $[[z]],$ $\cdot )$ are homomorphic, because we cannot
guarantee that $\Theta _{D_{N}}^{inv}$ and $D_{N}^{inv}[[z]]$ are bijective.
But, similar to the scalar-valued case, we can define a map $F$ from $\Theta
_{D_{N}}^{inv}$ to $D_{N}^{inv}[[z]]$ by

\strut

(4.5) $\ \ \ \ \ \ \ F\left( g(z)\right) =\frac{1}{z}\,\,g^{<-1>}(z),$ for
all $g\in \Theta _{D_{N}}^{inv}.$

\strut

Then, clearly, the map $F$ is injective. By replacing $\Bbb{C}$ to $D_{N}$
in the proof in [1], we can prove the above theorem.

\strut

\begin{definition}
Let $x\in (A,E)$ be a $D_{N}$-valued random variable such that $E(x)$ $\in $ 
$D_{N}^{-1}$, where $A$ $=$ $\times _{j=1}^{N}$ $A_{j}$ and $E$ $=$ $%
(\varphi _{1},$ ..., $\varphi _{N}).$ Then the S-transform of $x$ is defined
by

\strut 

$\ \ \ \ \ \ \ \ \ \ \ \ \ \ \ S_{x}(z)\overset{def}{=}\frac{1}{z}%
R_{x}^{<-1>}(z)$ \ in \ $D_{N}^{inv}[[z]].$
\end{definition}

\strut \strut \strut

By the previous theorem, we can get the following main result of this
chapter;

\strut

\begin{corollary}
Let $x$ and $y$ be $D_{N}$-valued random variables in our direct producted
noncommutative probability space $(A,$ $E)$ with amalgamation over the $N$%
-th diagonal algebra $D_{N},$ satisfying that $E(x),$ $E(y)$ $\in $ $%
D_{N}^{-1}.$ If $x$ and $y$ are free over $D_{N}$ in $(A,$ $E),$ then

\strut 

(4.6) \ \ \ \ \ \ \ \ $S_{xy}(z)=\left( S_{x}(z)\right) \left(
S_{y}(z)\right) ,$ \ in $D_{N}^{inv}[[z]].$
\end{corollary}

\strut

\begin{proof}
\strut Let $x$ and $y$ be $D_{N}$-valued random variables in $(A,$ $E)$ and
assume that $E(x)$ and $E(y)$ are in $D_{N}^{-1}.$ Then the R-transforms $%
R_{x}(z)$ and $R_{y}(z)$ of $x$ and $y$ have their composition inverse $%
R_{x}^{<-1>}(z)$ and $R_{y}^{<-1>}(z)$ in $D_{N}$ $[[z]].$ Also, we have
that $E(xy)$ $=$ $E(x)$ $E(y),$ by the $D_{N}$-freeness of $x$ and $y,$ and
hence $E(xy)$ is also contained in $D_{N}^{-1}$. So, $R_{xy}^{<-1>}(z)$
exists in $D_{N}$ $[[z]].$ Thus the S-transforms $S_{xy},$ $S_{x}$ and $%
S_{y} $ are well-defined in $D_{N}^{inv}$ $[[z]],$ with their $D_{N}$%
-constant terms $E(xy),$ $E(x)$ and $E(y),$ respectively. Observe that

\strut

$\ \ \ \ S_{xy}(z)=\frac{1}{z}R_{xy}^{<-1>}(z)=\frac{1}{z}\left( R_{x}\,\,\,%
\frame{*}_{D_{N}}\,\,\,R_{y}\right) ^{<-1>}(z)$

\strut

by the $D_{N}$-freeness of $x$ and $y$

\strut 

$\ \ \ \ \ \ \ \ \ \ \ =F\left( R_{x}\,\,\,\frame{*}_{D_{N}}\,\,\,R_{y}%
\right) =\left( F(R_{x})\right) \left( F(R_{y})\right) $

\strut

by the definition of the group homomorphism $F$ from $(\Theta
_{D_{N}}^{inv}, $ $\ \frame{*}_{D_{N}})$ to $(D_{N}^{inv}$ $[[z]],$ $\cdot )$

\strut

$\ \ \ \ \ \ \ \ \ \ \ =\left( \frac{1}{z}R_{x}^{<-1>}(z)\right) \left( 
\frac{1}{z}R_{y}^{<-1>}(z)\right) =\left( S_{x}(z)\right) \left(
S_{y}(z)\right) .$

\strut
\end{proof}

\strut

We finish this chapter by observing the following general case;\strut

\strut

\begin{remark}
In this remark, we will extend the above observation to a certain general
case. Let's assume that $D$ is a commutative unital algebra. Define a subset 
$\Theta _{D}$ of formal series $D[[z]]$ by

\strut 

$\ \ \ \ \ \ \ \Theta _{D}=\{\sum_{n=0}^{\infty }d_{n}z^{n}\in
D[[z]]:d_{0}=0_{D}\}.$

\strut 

Define the restricted $D$-valued boxed convolution \frame{*}$_{D}$ $:$ $%
\Theta _{D}$ $\times $ $\Theta _{D}$ $\rightarrow $ $\Theta _{D}$ by

\strut 

$\ \ \ \ \ \ \ \ \ \ \ \ \ \ \ \left( g_{1}\,\,\frame{*}_{D}\,\,g_{2}\right)
(z)=\sum_{n=1}^{\infty }d_{n}z^{n}$

with

$\ \ \ \ \ \ \ \ \ \ \ \ \ \ \ \ \ d_{n}=\underset{\pi \in NC(n)}{\sum }%
\left( d_{\pi }^{(1)}\right) \left( d_{Kr(\pi )}^{(2)}\right) ,$

where

$\ \ \ \ \ \ \ g_{1}(z)=\sum_{n=1}^{\infty }d_{n}^{(1)}\,z^{n}$ \ and \ $%
g_{2}(z)=\sum_{n=1}^{\infty }d_{n}^{(2)}z^{n}$

and

\strut 

$\ \ \ d_{\theta }^{(i)}=\underset{V\in \theta }{\Pi }d_{\left| V\right|
}^{(i)},$ \ for \ $i$ $=1,2,$ $\forall \theta \in NC(n),$ $\forall n\in \Bbb{%
N}.$

\strut 

Then, the algebraic structure $\left( \Theta _{D},\,\frame{*}_{D}\right) $
is a semigroup, because, for any $g_{1},$ $g_{2},$ $g_{3}$ $\in $ $\Theta
_{D},$ we can get $D$-valued random variables $x_{1},$ $x_{2}$ and $x_{3}$
in the $D$\textbf{-free} noncommutative probability space $(A_{1},$ $E_{1}),$
$(A_{2},$ $E_{2})$ and $(A_{3},$ $E_{3})$ with amalgamation over $D,$
respectively, such that $D$ commutes with $A_{1},$ $A_{2}$ and $A_{3},$
respectively, and $g_{j}(z)$ $=$ $R_{x_{j}}(z),$ for $j$ $=$ $1,$ $2,$ $3,$
where $R_{x_{j}}(z)$ is defined by

\strut 

$\ \ \ \ \ \ \ \ \ \ \ R_{x_{j}}(z)=\sum_{n=1}^{\infty }\left(
k_{n}^{E}(x_{j},...,x_{j})\right) \,z^{n}$ \ in \ $\Theta _{D},$

\strut 

for $i=1,2,3,$ where $E=E_{1}*E_{2}*E_{3}$ on $A_{1}*_{D}A_{2}*_{D}A_{3}.$
Indeed,

\strut 

$\ \ \ \left( g_{1}\,\,\frame{*}_{D}\,\,g_{2}\right) \,\frame{*}%
_{D}\,\,g_{3}=R_{x_{1}x_{2}x_{3}}=g_{1}\,\,\frame{*}_{D}\,\,\left( g_{2}\,\,%
\frame{*}_{D}\,\,g_{3}\right) .$

\strut 

So, $\left( \Theta _{D},\,\,\frame{*}_{D}\right) $ is a semigroup. Define a
subset $\Theta _{D}^{inv}$ in $\Theta _{D}$ by

\strut 

$\ \ \ \ \ \ \ \ \ \ \ \Theta _{D}^{inv}=\{\sum_{n=1}^{\infty
}d_{n}\,z^{n}\in \Theta _{D}:d_{1}\in D_{inv}\},$

\strut 

where $D_{inv}=\{d\in D:d$ is invertible in $D\}.$ Then $\left( \Theta
_{D}^{inv},\,\,\frame{*}_{D}\right) $ is a group with its \frame{*}$_{D}$%
-identity $1_{D}\cdot z.$ If we define

\strut 

$\ \ \ \ \ \ \ \ \ D[[z]]_{inv}=\{g\in D[[z]]:g$ is ($\cdot $)-invertible$\},
$

\strut 

where ($\cdot $) is the usual multiplication on $D[[z]].$ Then $\left(
D[[z]]_{inv},\cdot \right) $ is a group, too. By defining

\strut 

$\ \ \ \ \ \ \ \ \ \ \ \ \ \ \ \ \ \ \ \ \ F_{D}:\Theta
_{D}^{inv}\rightarrow D[[z]]_{inv}$

by

$\ \ \ \ \ \ \ \ \ \ \ \ \ \ \ \ \ \ \ \ \ \ \ F_{D}(g)=\frac{1}{z}%
g^{<-1>}(z),$

\strut 

where $g^{<-1>}$ is the composition inverse of $g$ in $D[[z]],$ we can
verify that $\left( \Theta _{D}^{inv},\,\frame{*}_{D}\right) $ and $\left(
D[[z]]_{inv},\,\cdot \right) $ are homomorphic. So, we can define the
S-transform $S_{x}(z)$ of a $D$-valued random variable $x$ in the
noncommutative probability space $(B,$ $E_{D})$ with amalgamation over $D,$
with $E_{D}(x)$ $\in $ $D_{inv},$ where $D$ commutes with $B,$ by

\strut 

$\ \ \ \ \ \ \ \ \ \ \ \ \ \ \ \ \ \ \ \ \ \ \ \ \ S_{x}(z)=F\left(
R_{x}(z)\right) ,$

\strut 

where $R_{x}(z)$ is the R-transform of $x\in (B,E_{D})$ in $\Theta
_{D}^{inv}.$ And if $x$ and $y$ are $D$-free random variables in $(B,$ $%
E_{D})$ such that $E_{D}(x)$ and $E_{D}(y)$ are contained in $D_{inv},$ then

\strut 

$\ \ \ \ \ \ \ \ \ \ \ \ \ \ S_{xy}(z)=\left( S_{x}(z)\right) \left(
S_{y}(z)\right) $ \ \ in \ \ $D[[z]]_{inv}.$
\end{remark}

\strut \strut \strut 

\strut

\strut \strut

\strut \textbf{References}

\bigskip

\strut

\begin{quote}
{\small [1]\ \ \ A. Nica, R-transform in Free Probability, IHP course note.}

{\small [2]\ \ \ A. Nica, R-transforms of Free Joint Distributions and
Non-crossing Partitions, J. of Func. Anal, 135 (1996), 271-296.\strut }

{\small [3]\ \ \ D.Voiculescu, K. Dykemma and A. Nica, Free Random
Variables, CRM Monograph Series Vol 1 (1992).\strut }

{\small [4] \ F. Radulescu, Singularity of the Radial Subalgebra of }$%
L(F_{N})${\small \ and the Puk\'{a}nszky Invariant, Pacific J. of Math, vol.
151, No 2 (1991)\strut , 297-306.\strut \strut }

{\small [5] \ \ I. Cho, Toeplitz Noncommutative Probability Spaces over
Toeplitz Matricial Algebras}$,${\small \ (2002), Preprint.\strut }

{\small [6] \ \ I. Cho, The Moment Series of the Generating Operator of }$%
L(F_{2})*_{L(K)}L(F_{2})${\small , (2003), Preprint.}

{\small [7] \ \ I. Cho, Random Variables in a Graph }$W^{*}${\small %
-Probability Space, (2003), Preprint. }

{\small [8]\ \ \ I. Cho, The Moments of Block Operators of a Group von
Neumann Algebra, (2005), Preprint.}

{\small [9] \ \ K. J. Horadam, The Word Problem and Related Results for
Graph Product Groups, Proc. AMS, vol. 82, No 2, (1981) 157-164.}

{\small [10]\ R. Speicher, Combinatorics of Free Probability Theory IHP
course note }

{\small [11] R. Speicher, Combinatorial Theory of the Free Product with
Amalgamation and Operator-Valued Free Probability Theory, AMS Mem, Vol 132 ,
Num 627 , (1998).}
\end{quote}

\strut

\end{document}